\documentclass[trsc,nonblindrev]{informs3} 
\graphicspath{{../}}
\OneAndAHalfSpacedXI 

\usepackage{amsmath,verbatim,amssymb,amsfonts,amscd,mathtools,xcolor,latexsym}
\usepackage{algorithm}
\usepackage{algcompatible} 
\usepackage[noend]{algpseudocode}
\usepackage[stable]{footmisc}
\usepackage[T1]{fontenc}
\usepackage{comment}
\usepackage{mathrsfs}
\usepackage{wrapfig}
\usepackage{caption} 
\usepackage{subfigure}
\usepackage{lscape}
\usepackage{graphicx,tabularx,tabulary}
\usepackage{booktabs}
\usepackage[breaklinks=true, colorlinks=true, linkcolor={blue!90!black}, citecolor={blue!90!black}, urlcolor={blue!90!black}]{hyperref}
\usepackage{cleveref} 
\usepackage{thmtools} 
\usepackage{thm-restate} 
\usepackage{mathtools}

\usepackage{tikz-cd}
\usepackage{mhchem} 

\newcommand{\E}{\mathbb{E}}

\newcommand{\cL}{\mathcal{L}}

\newcommand{\cE}{\mathcal{E}}

\newcommand{\cG}{\mathcal{G}}
\newcommand{\cN}{\mathcal{N}}

\newcommand{\cR}{\mathcal{R}}

\newcommand{\cS}{\mathcal{S}}
\newcommand{\cU}{\mathcal{U}}

\newcommand{\cV}{\mathcal{V}}
\newcommand{\bZ}{\mathbb{Z}}
\newcommand{\cF}{\mathcal{F}}

\newcommand{\rtw}{\text{RTW}}

\newcommand{\obj}{\text{obj}}

\newcommand{\wt}[1]{\widetilde{#1}}

\newcommand{\ol}{\overline}
\newcommand{\ul}{\underline}
\newcommand{\mtc}{\mathcal}

\newcommand{\rtau}{{\hat{\tau}^{\textrm{rel}}}}

\newcommand{\Tdiff}{T_{\textrm{diff}}}
\newcommand{\Gdiff}{G_{\textrm{diff}}}
\newcommand{\Gmax}{G_{\textrm{max}}}
\newcommand{\Gmin}{G_{\textrm{min}}}

\newcommand{\gamfull}{\gamma_{\text{full}}}
\newcommand{\gamext}{\gamma_{\text{ext}}}

\newcommand{\mL}{\mathcal{L}}

\newcommand{\rT}{\wt{T}}

\newcommand{\ut}{\underline{t}}
\newcommand{\ot}{\overline{t}}

\usepackage[normalem]{ulem}
 
\newcommand{\beq}{\begin{equation}}
\newcommand{\eeq}{\end{equation}}
\newcommand{\bdm}{\begin{displaymath}}
\newcommand{\edm}{\end{displaymath}}
\newcommand{\ba}{\begin{aligned}}
\newcommand{\ea}{\end{aligned}}

\theoremstyle{EX}

\newtheorem{observation}{Observation}[section]

\newcommand{\assref}[1]{Assumption~\ref{#1}}

\newcommand{\propref}[1]{Proposition~\ref{#1}}
\newcommand{\obsref}[1]{Observation~\ref{#1}}
\newcommand{\defref}[1]{Definition~\ref{#1}}
\newcommand{\secref}[1]{Section~\ref{#1}}
\newcommand{\tabref}[1]{Table~\ref{#1}}
\newcommand{\figref}[1]{Figure~\ref{#1}}
\newcommand{\algoref}[1]{Algorithm~\ref{#1}}
\newcommand{\linref}[1]{Line~\ref{#1}}
\newcommand{\expref}[1]{Example~\ref{#1}}
\newcommand{\DT}{(DT)}

\DeclarePairedDelimiterX\Set[2]{\lbrace}{\rbrace}%
 { #1 \colon\mathopen{} #2 }

\usepackage{enumitem}

\usepackage{endnotes}
\let\footnote=\endnote

\usepackage{natbib}
 \bibpunct[, ]{(}{)}{,}{a}{}{,}%
 \def\bibfont{\small}%
\TheoremsNumberedThrough     
\ECRepeatTheorems

\EquationsNumberedThrough    

\begin{document}


\RUNTITLE{Vehicle Platooning on Tree Networks}

\TITLE{Coordinated Vehicle Platooning on Tree Networks:
Efficient Time Discretization and Strengthened Formulation}

\ARTICLEAUTHORS{%
\AUTHOR{Fengqiao Luo}
\AFF{Uber Technologies, Inc.}
\AFF{Northwestern University, Evanston, IL, 60208,
\EMAIL{fengqiaoluo2014@u.northwestern.edu}}
} 

\ABSTRACT{We consider the coordinated vehicle platooning problem on a tree network
with time constraints while the routes of vehicles are given. 
The problem is to coordinate the
departure time of each vehicle to enable platoon formation hence maximizing the total fuel saving.
For this problem setting, relative time windows
can be defined for all vehicles to which an efficient time discretization can be applied.
This property leads to a tight mixed-integer linear program reformulation
as compared to the continuous-time formulation involving big-M coefficients proposed in our previous work. 
It is demonstrated by systematic numerical experiments that the reformulation outperforms
the continuous-time formulation for this family of problem instances.
}


\maketitle

\section{Introduction}
Vehicle platooning is a navigation and control technology under field test that is proposed to be
implemented in an intelligent transportation system (ITS), especially with applications to trucks to improve
safety, increase energy efficiency, reduce cost and CO$_2$ emission, and increase road-usage efficiency. 
In a platoon of trucks on the road, multiple trucks drive in tandem without being physically connected,
and one truck follows another with a small inter-vehicle distance, e.g. 6m$\sim$10m. 
Vehicle (truck) platooning is often associated with autonomous driving, 
where only the lead vehicle needs to be operated by a driver and the trailing vehicles are 
automatically (or semi-automatically) controlled. 
Vehicle platooning technology includes automatic braking and speed control, 
and the system is able to react much faster than a human, improving safety and reducing the 
likelihood of collisions. As a consequence, a much shorter inter-vehicle distance is needed
to guarantee safety compared to the human-driving case, hence increasing the road-usage
efficiency and throughput. Trucks (especially trailing trucks) in a platoon can drive at lower 
fuel (energy) cost compared to solo driving due to reduced aerodynamic drag incurred. 
Extra energy saving comes from adaptive cruise control.  

According to the research \cite{muratori2017potentials} conducted by the National Renewable 
Energy Laboratory based on data collected from in-car devices, 
65\% of the miles travelled by trucks could be platooned, 
resulting in a 4\% reduction in total truck fuel consumption. 
The fuel saving factor is impacted by the inter-vehicle distance, 
the position of the vehicle within the platoon and the travel speed. 
Results demonstrated a wide range of fuel savings \cite{McAuliffe2018influence} 
 -- with the lead vehicle saving up to 10\% 
at the closest separation distances, the middle vehicle saving up to 17\%, 
and the trailing vehicle saving up to 13\%.
On average the lead truck demonstrated fuel saving up to 5.3\% and the trailing truck saved up to 9.7\%.
While platooning improved fuel economy at all speeds, 
travel at 55 mph resulted in the best overall miles per gallon, referring to \cite{2014-eff-platn-fuel-speed-dist-mass,McAuliffe2018influence}. 
Innovative hardwares and technologies are needed to enable platoon formation, 
splitting and maintenance with safety.
According to Peloton Technology \cite{pelotonTech}, 
a pioneer company in the industry of connected vehicle technology,
in-car devices for platooning consists of the following five major components: platooning control unit, 
wireless communication system, radar-based collision mitigation, driver controls
and system display. 

The planning for truck platooning at scale over a road network can be divided into three levels 
by \cite{hoef2018thesis-coord-veh-platn}: 
the \textit{strategic level}, the \textit{tactical level} and the \textit{operational level}.
At the strategic layer, the strategies have been generated by a central system named platoon coordinator
offered by a platoon service provider, an organization that provides crucial platooning services shared 
between road transport providers for services such as certification, insurance and coordination, referring to
\cite{janseen2015truck-platn}.
The platoon coordinator can be a cloud that receives data from traveling or static vehicles
spread over the network via base stations, and then compute routes, departure time and speed profiles
for vehicles to form platoon when they merge on a shared road segment. The computation
should utilize the forecast of traffic and its impact on the travel time. 
The information out from the computation will be sent to vehicles via base stations.
At the tactical layer, the platoon manager (usually played by the lead vehicle) controls the formation, 
splitting and reorganization of the platoon based on the plans provided by the strategic layer.
This can be achieved using V2V communication. At the operational layer, the in-car devices
from each member vehicle in a platoon start controlling the vehicle with the reference speed 
and inter-vehicle distance set by the platoon manager.

It is worth to note that truck platooning is more than an experimental technology. 
Industrialization of the technology has made steady progress, 
and the first generation of it is close for commercial launch.
The truck platooning market was valued at USD 1405 million in 2020, 
and it is expected to reach USD 6077 million by 2026, registering a CAGR of about 32\% 
during the forecast period (2021-2026) according to \cite{platnMarketGrow2021}.
North America is the current largest market and 
Europe is the fastest growing market in which the technology development is led by Sweden
and Germany. Major players are Peloton Technology, Daimler AG, Volvo Group, 
DAF Trucks (Paccar Company) and Scania AB, etc.

\subsection{Literature review}
The paper contributes to the development of strategic solutions for vehicle platooning
over the road network, with focus on high-level time coordination to enable platooning.
We review related literature on the development of platooning technology 
and the research on the strategic level. 

\subsubsection*{Platooning technology development}
Implementation of truck platooning technology depends on the development
of hardware and software used for adaptive cruise control and wireless communication.
We refer to \cite{2010-operating-platoons}, \cite{tsugawa2013overview} and \cite{Tsugawa2016} 
as an overview of platooning technology and observed energy saving, and 
refer to \cite{2012HobertMaster} and \cite{2012-coord-veh-plat-control} for the impact of 
wireless communication protocol on platoon formation. 
Speed control algorithms applied during the process of platoon formation 
play an essential role in achieving fuel saving. The control theory and algorithms
for maximizing the fuel saving at the vehicle level are investigated in 
\cite{shida2009development}, \cite{2013-fuel-eff-platn-HDV}, \cite{Hoef2015:itsc:15}, \cite{Zhang2019}, \cite{Gong2016} and \cite{Ye2019}.
The control algorithms for maintaining a platoon at energy efficiency is studied in \cite{2011-truck-platoon-control}.
A variety of experiments and road tests have been conducted to gain a comprehensive understanding 
of the amount of fuel saving in different situations and the impact factors (e.g., speed, gap, wind, and slope, etc.),
referring to \cite{davila2013sartre}, \cite{alam2010experimental}, \cite{browand2004fuel}, \cite{2003-truck-eng-opt},
\cite{2000-fuel-reduction-platoon}, \cite{2004-fuel-saving-two-trucks} and \cite{2012-cost-curves-HDV}.
The data from experiments is used to design simulators and control algorithms.
At the network level, \cite{2014-fuel-sav-platn-HDV-posit-data} analyzed sparse data of truck position
to identify potential platooning opportunities. Besides of fuel saving, there are studies concerning
the impact of truck platooning on the highway throughput, referring to
\cite{vander2002effects}, \cite{shladover2012impacts} and \cite{2012-plat-mitig-delay-improve-traff-flow},
and on the roadside hardware and safety, referring to \cite{roadside2019}.

\subsubsection*{Strategic solution for mining platooning opportunities}
On the strategic level, the central controller (either a control tower or a cloud) needs to create 
a plan for route design and travel schedule coordination of multiple vehicles to optimize the overall
fuel consumption (of the system) by enabling achievable platooning. Note that the scale of region managed by the 
central controller can vary. In particular, the region can cover a few highway junctions, referring to 
\cite{xiong2021-dyn-prog-veh-platn}, or it can be at the level of a road network that distributes over multiple states,
referring to \cite{luo2018-veh-platn-mult-speeds} and \cite{luo2020-route-then-sched-vehicle-platn}. 
In general, there could be central controllers managing regions at different scale that form 
a top-down hierarchy. There is a series of investigation on optimal plan for platooning of
a handful of vehicles using mathematical programming based models,
referring to \cite{Baskar2009a}, \cite{2016-coord-platn-routing}, \cite{Hoef2015a:acc:15}. 
Some models integrate route design, time coordination and speed control to find an optimal
plan to minimize the overall fuel consumption
during the entire process, referring to \cite{Kammer2013}, \cite{Liang2014}, \cite{Hoef2015a:acc:15} and \cite{Baskar2013}. 
The quantitative impact of traffic flow and road capacity on the platoon
formation is also incorporated in some studies, referring to \cite{Baskar2009:itsc:09}, and \cite{Baskar2013}.
These models often contain nonlinear terms
in the objective to count the impact of speed change on the fuel saving factor.
Some other models focus on high-level planning using mixed-integer linear programs 
with an overall fuel saving factor that only discriminates lead and trailing vehicles but
not intending to model the platoon formation process, and hence these models are better
at scale, referring to \cite{Larson2014d}, \cite{Larsson2015}, \cite{luo2018-veh-platn-mult-speeds} 
and \cite{luo2020-route-then-sched-vehicle-platn}.
There are some other studies using local search to propose routes suitable for platooning, 
referring to \cite{Doremalen2014}, \cite{Larson2013b:itsc:13}. 
As the platooning process involves speed change
and catch-up, dynamic programming is also a useful tool to model it at the tactical level, referring to 
studies in \cite{Valdes2012} and \cite{xiong2021-dyn-prog-veh-platn} as examples. 
Recently, \cite{xiong2021-dyn-prog-veh-platn} developed a novel stochastic dynamic programming model
to optimize the platooning process focusing on the case when multiple vehicles merge at the highway junction. 
\cite{luo2020-route-then-sched-vehicle-platn} established a novel repeatedly route-then-schedule framework 
for solving coordinated vehicle platooning problem with time constraints for hundreds of vehicles 
distributed over the road network. 

\subsection{Problem statement}
In the tree-based fixed-route coordinated vehicle platooning (tree-FRCVP) problem, 
we are given a set $\cV$ of vehicles. For each vehicle $v\in\cV$, 
the following information is given: the origin node $O_v$, destination node $D_v$, 
earliest possible departure time (from the origin) $T^O_v$, and the deadline of arriving 
time (at the destination) $T^D_v$. In addition, the route $\cR_v$ of
each vehicle $v$ is given, which consists of arcs (directed edges) that form a path from
the origin to the destination of $v$. Note that vehicles can have shared
arcs in their routes. Let $\cG=\cup_v\cR_v$ be the directed graph of vehicle routes (the route graph, in short). 
The objective is to find an optimal departure time 
of each vehicle to maximize fuel saving of the vehicle system 
by forming platoons on arcs along the vehicle routes. Multiple 
vehicles can form a platoon on an arc if they enter the starting node of an arc at the same time. 
We focus on the deterministic problem setting in the paper. 
Incorporation of the uncertainty in the travel time and 
and the impact of traffic requires further development, which is briefly discussed in
Section~\ref{sec:discussion}. 
The following assumption summarizes the problem setting.
Note that the route graph $\cG$ can be an arbitrary directed graph in general.
However, the focus of this work is on the case that $\cG$ is a tree-type network.
\begin{assumption}\label{ass:polytree}
The directed graph $\cG=\cup_{v\in\cV}\cR_v$ is a polytree, i.e., the undirected counterpart of $\cG$ is a tree.
\end{assumption}
An example of a polytree network is given in \figref{fig:poly-path}.
We investigate the coordinated vehicle platooning (CVP) problem based on this assumption as we 
have discovered that, there are interesting properties of the CVP problem on a tree-network that can
lead to a tight MILP formulation which outperforms the formulation developed in our previous work on solving
CVP problems in this particular case. Even for practical applications, this assumption can be well rationalized: 
The platooning technology is usually applied at major high-way segments in a region. 
Multiple vehicles join a same high-way segment simultaneously and form a platoon until
they leave each other. It is unlikely or impractical for two vehicles to join a platoon
on separate high-way segments while leave each other in the middle, which induces cycles
in $\cG$. Therefore, the special case of $\cG$ being a polytree can represent an important family of CVP problem 
instances.

\begin{assumption}\label{ass:no-pause}
\emph{(1)} For any given arc in $\cG$, the speed on the arc is a given parameter, which is the same for all vehicles
when traveling on this arc, but it can be arc dependent.
\emph{(2)} The fuel saving factor is a constant, and it is $\sigma^l$ for a lead vehicle and $\sigma^f$ for a trailing vehicle. 
\end{assumption}
Throughout the paper, we also make the following assumption to regularize the tree-FRCVP instance
and the route of each vehicle:
\begin{assumption}\label{ass:irr}
(1) The tree-FRCVP problem instance is irreducible, i.e., 
there does not exist a partition $(\cV_1,\cV_2)$ of $\cV$ such that 
$(\cup_{v\in\cV_1}\cR_v)\cap(\cup_{v\in\cV_2}\cR_v)=\emptyset$.
\end{assumption}
Note that if a tree-FRCVP instance is reducible, it can be
decomposed into two independent instances and hence can be studied separately.

\subsection{Our contributions}
In this paper, we focus on the time coordination problem (scheduling problem) of 
coordinated vehicle platooning over a tree (road) network when the vehicle routes are given 
which is an algorithmic component at the strategic layer. Our contribution can be summarized as follows:
For tree-network instances,
we have discovered that there exists an efficient discretization of time windows.
In this case, a relative time window (RTW)
can be defined for each vehicle, and the time discretization can be naturally performed according to
the intersections among all the RTW's. The time coordination problem is then 
equivalent to the task of assigning vehicles to feasible time buckets, which can be formulated as  
a tight mixed-integer linear program with less variables and with(out) big-M coefficients
as compared to the continuous-time formulation proposed in our previous work.
We have developed better formulations for two cases depending on whether the 
stop-and-wait option is permitted at intermediate nodes along the route of a vehicle.
Numerical investigation shows the new formulation leads to better performance
compared to the continuous time formulation.

\subsection{Connection to the previous work}
The fixed-route coordinated vehicle platooning (FRCVP) problem 
in which the route graph $\cG$ is not necessarily a polytree,
is closely connected to the general coordinated vehicle platooning (GCVP) problem. 
In the GCVP problem, it requires to determine an optimal route and departure time for every vehicle to 
minimize the total fuel cost, whereas in the FRCVP problem vehicle routes are given. 
We have developed a repeated route-then-schedule
algorithmic framework to decompose the GCVP problem into routing sub-problems
and scheduling sub-problems in a heuristic approach, referring to \cite{luo2020-route-then-sched-vehicle-platn},
which overcomes the high complexity of the integrated formulation developed in \cite{luo2018-veh-platn-mult-speeds}
that has a large number of variables and big-M coefficients for medium and large problem instances.
The FRCVP problem is the scheduling sub-problem in this framework. Furthermore, the FRCVP problem
worths an independent investigation, as vehicle routes are sometimes pre-determined (e.g., adapting the shortest path) 
in practice. The tree-FRCVP problem is a special case of the FRCVP problem with $\cG$ being a polytree. 
 
\section{A preprocess of problem instances}
\label{sec:categ-FRCVP}
\begin{figure}
\centering
\subfigure[]{\includegraphics[width=0.49\linewidth,trim=3.5cm 16cm 4cm 3cm]{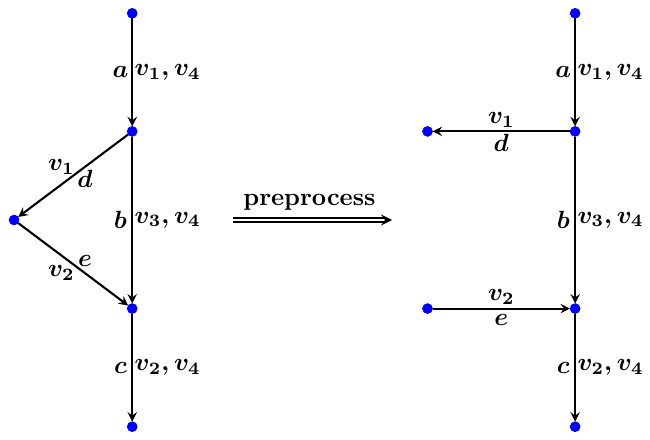}}
\hfill
\subfigure[]{\includegraphics[width=0.49\linewidth,trim=3.5cm 16cm 4cm 3cm]{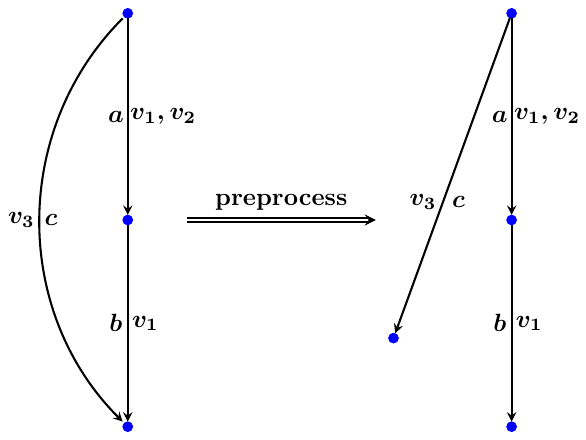}}
\captionsetup{justification=centering}
\caption{\footnotesize The network before and after the node-extension preprocess. 
In the above figures, if an arc $e$ is mapped to $e^\prime$, they are labeled by a same letter.
The two examples demonstrate that some FRCVP instances having (directed or poly-) cycles in the 
original digraph can be transformed into polytree networks using \algoref{alg:node-ext}.}
\label{fig:preprocess}
\end{figure}

There are some route graphs which originally contain (directed or poly-) cycles\footnotemark[1]
but they can be eliminated after a node-extension preprocess given in \algoref{alg:node-ext}. 
The idea of \algoref{alg:node-ext} is to first create a separate copy $\cR^\prime_v$ of the route $\cR_v$
for each $v\in\cV$. If an arc $e$ in the original route graph $\cG$ is on the routes of
a set $\cV_e$ of vehicles. The algorithm will create $|\cV_e|$ copies of $e$,
one in each $\cR^\prime_v$ for $v\in\cV_e$. Then the algorithm merges these copies
into a single arc. After this preprocess, it could happen that some cycles in the original route graph $\cG$
are broken in the processed route graph $\cG^\prime$. 
We now use the route graph $\cG$ in \figref{fig:preprocess}(a) to illustrate the algorithm.
In the loop at \linref{lin:preproc-loop1} of the algorithm, 
it creates the following routes for vehicles $\{v_1,v_2,v_3,v_4\}$:
\bdm
\cR^\prime_{v_1}=\{a_1,d_1\},\quad\cR^\prime_{v_2}=\{e_2,c_2\}, \quad
\cR^\prime_{v_3}=\{b_3\},\quad\cR^\prime_{v_4}=\{a_4,b_4,c_4\},
\edm 
where $\{a_1,a_4\}$ are a copies of $a$, $\{b_3,b_4\}$ are copies of $b$,
$\{c_2,c_4\}$ are copies of $c$, $d_1$ is a copy of $d$ 
and $e_2$ is a copy of $e$ from $\cG$, respectively. 
So $\varphi(a_1)=\varphi(a_4)=a$ and $\varphi(b_3)=\varphi(b_4)=b$, etc.
In the loop at \linref{lin:preproc-loop2} of the algorithm,
it merges $\{a_1,a_4\}$ to a single arc $a^\prime$, 
$\{b_3,b_4\}$ to a single arc $b^\prime$,
and $\{c_2,c_4\}$ to a single arc $c^\prime$. 
The processed figure $\cG^\prime$ consists of arcs 
$\{a^\prime, b^\prime, c^\prime, d_1, e_2\}$,
which are relabeled to be $\{a, b, c, d, e\}$
to match with the arc labels in $\cG$. Note that
the poly-cycle $\{d,e,b\}$ in $\cG$ is broken in $\cG^\prime$.

\footnotetext[1]{A directed cycle consists of arcs that are in the same direction and 
a poly-cycle consists of arcs that are in different directions.} 

\begin{observation}\label{obs:map-phi}
Let $\cG$ and $\cG^\prime$ be the input and output digraphs in \algoref{alg:node-ext}, respectively,
and let $\cE$ and $\cE^\prime$ be their sets of arcs.
There exists a bijective mapping $\phi$ from $\cE$ to $\cE^\prime$ such that for every $v\in\cV$
the arcs of $\cR_v$ (in $\cG$) are mapped into a route $\cR^\prime_v$ in $\cG^\prime$ linked in
the same order. While the number of nodes in $\cG^\prime$ is no less than that in $\cG$.
\end{observation}
\proof{Proof.}
Consider the mapping $\varphi$ constructed in the algorithm. For every $e\in\cE$, 
if $e$ is in $\cR_v$ for some $v$, the algorithm will create a new copy $\hat{e}^v$ of $e$
in $\cR^\prime_v$ and set $\varphi(\hat{e}^v)=e$. This indicates that $\varphi^{-1}(e)$ is the set of
copies of $e$ in $\cR^\prime_v$ for each $v\in\cV_e$. By the algorithm, all elements in $\varphi^{-1}(e)$
are merged into a single arc $e^\prime$ in $\cE^\prime$. So the bijective mapping $\phi$
satisfies $\phi(e)=e^\prime$, and the order of arcs in a route is preserved by $\phi$. 
\halmos
\endproof

\begin{algorithm}
{\footnotesize 
\caption{\footnotesize Node-extension preprocess.}\label{alg:node-ext}
\begin{algorithmic}[1]
\State{\bf Input}: the vehicle set $\cV$ and directed graph $\cG=\cup_{v\in\cV}\cR_v$.
\State{\bf Output}: a transformed directed graph $\cG^\prime$.
\State{Initialize a mapping $\varphi$, with the mapping rule specified later.}
\For{$v\in\cV$} \label{lin:preproc-loop1}
	\State{Let $k=|\cR_v|$.}
	\State{Create a new directed path $\cR^\prime_v=(\hat{e}^v_1,\ldots,\hat{e}^v_k)$.}
	\State{Set $\varphi(\hat{e}^v_i)=e^v_i$ for $i=1,\ldots,k$, where $\{e^v_i:i=1,\ldots,k\}$
	 are arcs along $\cR_v$, i.e., $\cR_v=(e^v_1,\ldots,e^v_k)$.}
	\State{Create a new set of distinct nodes $(s^v_1,\ldots,s^v_{k+1})$ such that $\hat{e}^v_{i}=(s^v_i,s^v_{i+1})$ for $i=1,\ldots,k$.}
\EndFor
\For{$e\in\cG$}  \label{lin:preproc-loop2}
	\State{Merge the elements in $\varphi^{-1}(e)$ to be a single arc.}
	\State{Merge corresponding end nodes of directed edges in $\varphi^{-1}(e)$.}
\EndFor
\State{After finishing the above loop, the routes $\cR^\prime_v\;\forall{v\in\cV}$ have been glued up
by the shared directed edges that are mapped to a same element in $\cG$. Let the resulting directed graph be $\cG^\prime$.}
\State{\Return{$\cG^\prime$}.}
\end{algorithmic}
}
\end{algorithm}

\section{A Vehicle-to-time-bucket Assignment Formulation}\label{sec:assign-form-tree-FRCVP}
The continuous-time formulation of FRCVP established in \cite{luo2020-route-then-sched-vehicle-platn}
uses big-M coefficients. (The formulation is given in the appendix for completeness.)
In this section, we establish a big-M free formulation based on an efficient time-discretization technique
for tree-FRCVP instances which may improve the computational performance of solving the tree-FRCVP. 
Notations used to develop this formulation are summarized in Table~\ref{tab:param1}. 
In this section, we consider the case that all vehicles are not permitted to stop at intermediate nodes
to wait for other vehicles to form a platoon. They can only coordinate their departure times
at their origin nodes to enable platooning later on. 
\begin{assumption}\label{ass:no-pause}
Once a vehicle has departed from the origin, it is not allowed to stop and wait at any intermediate node along the route. 
\end{assumption}
The above assumption is rational for most practical cases since: 1. The energy cost of restarting
a truck can overwhelm the amount saved by platooning. 2. Building waiting areas along
the highway to enable the stop-and-wait option is impractical. 3. The stop-and-wait
option can increase the risk of having traffic accidents. 
Note that this assumption is also used in \cite{luo2020-route-then-sched-vehicle-platn}.
In \secref{sec:allow-waiting}, we extend our formulation development to cover the case when this 
assumption is removed.

The following definitions are useful
in describing the time-discretization technique which is the basis for our formulation. 
\begin{definition}[dir-path shared by vehicles]\label{def:share-path}
Let $P$ be a dir-path (directed path) on the graph $\cG$, and $\cV^\prime$ be a subset of vehicles. 
Vehicles in $\cV^\prime$ \emph{share} the dir-path $P$ if $P\subseteq\cR_v
\;\forall v\in\cV^\prime$. We use the notation $\cV(P)$ to denote the set 
of vehicles that share the dir-path $P$.
\end{definition}

\begin{definition}[time window]\label{def:time-wind}
Denote $\underline{t}_{u,i}$ and $\overline{t}_{u,i}$ ($\forall u\in\cV,\forall i\in\cN_u$)
as the lower and upper bound of the arrival time at node $i$ for vehicle $u$, where 
$\underline{t}_{u,i}=T_u^{O}+\sum_{(r,s)\in\cR_u(O_u,i)}T_{rs}$
and $\overline{t}_{u,i}=T_u^{D}-\sum_{(r,s)\in\cR_u(i,D_u)}T_{rs}$.
For any vehicle $v\in \cV$ and any node $i\in\cN_v$, the \emph{time window}
of the vehicle $v$ \emph{at} node $i$ is defined as $[\underline{t}_{v,i},\overline{t}_{v,i}]$.
\end{definition}

\begin{table}
  \footnotesize
  \caption{Sets, the graph representation, and parameters. \label{tab:param1}}
  \begin{tabularx}{\textwidth}{lX}
    \toprule
    Set and Graph & Definition\\
    \midrule
			$\cV$ & the set of vehicles \\
			$\cR_v$ & the route (a set of arcs that form a dir-path from $O_v$ to $D_v$) of vehicle $v$, for $v\in\cV$ \\
			$\cR$ & the union of all vehicle routes, $\cR=\cup_{v\in\cV}\cR_v$ \\
			$e$ & the short notation of an arc, e.g., $e=(i,j)$ \\
			$\cR_v(i,j)$ & the sub-route of vehicle $v$ from the node $i$ to the node $j$ along $\cR_v$ 
					(the two nodes can be non-consecutive but $i$ must be ahead of $j$) \\
			$\cV_e$ & the set of vehicles that share the arc $e$ along their routes \\
			$\cN_v$ & the set of nodes in $\cR_v$ \\
			$\cN$ & the union of all nodes, i.e., $\cN=\cup_{v\in\cV}\cN_v$ \\ 
			$\cE$ & the set of all arcs, i.e., $\cE=\cup_{v\in\cV}\cR_v$ \\
			$\cG$ & the directed graph defined by the set $\cN$ of nodes and the set $\cE$ of arcs \\
			$\cV(P)$ & the subset of vehicles that share the dir-path $P$ (see Definition~\ref{def:share-path}) \\
    \midrule\\[-.35cm]
    Parameter & Definition\\
    \midrule
			$C_{i,j}$ & fuel cost of traversing arc $(i,j)\in \cE$ for each vehicle\\
			$\sigma^l$, $\sigma^f$ & fuel saving rate for the lead and following vehicle in a platoon\\
			$O_v, D_v$ & origin and destination nodes for vehicle $v\in\cV$  \\
			$T_v^O$ & earliest allowable departure time of vehicle $v$ from its origin\\
			$T_v^D$ & latest allowed arrival time of vehicle $v$ at its destination\\
      			$T_{i,j}$ & time cost of traversing arc $(i,j)\in \cE$ for all vehicles\\
			$\lambda$ & the maximum number of vehicles that are allowed to form a single platoon on each road link \\
			$\underline{t}_{u,i}$ and $\overline{t}_{u,i}$ & given in Definition~\ref{def:time-wind} for $u\in\cV$ and $i\in\cR_u$ \\
    \bottomrule
  \end{tabularx}
\end{table}

\begin{figure}
\centering
    \includegraphics[width=0.5\linewidth,trim=0.2cm 13cm 0cm 4cm]{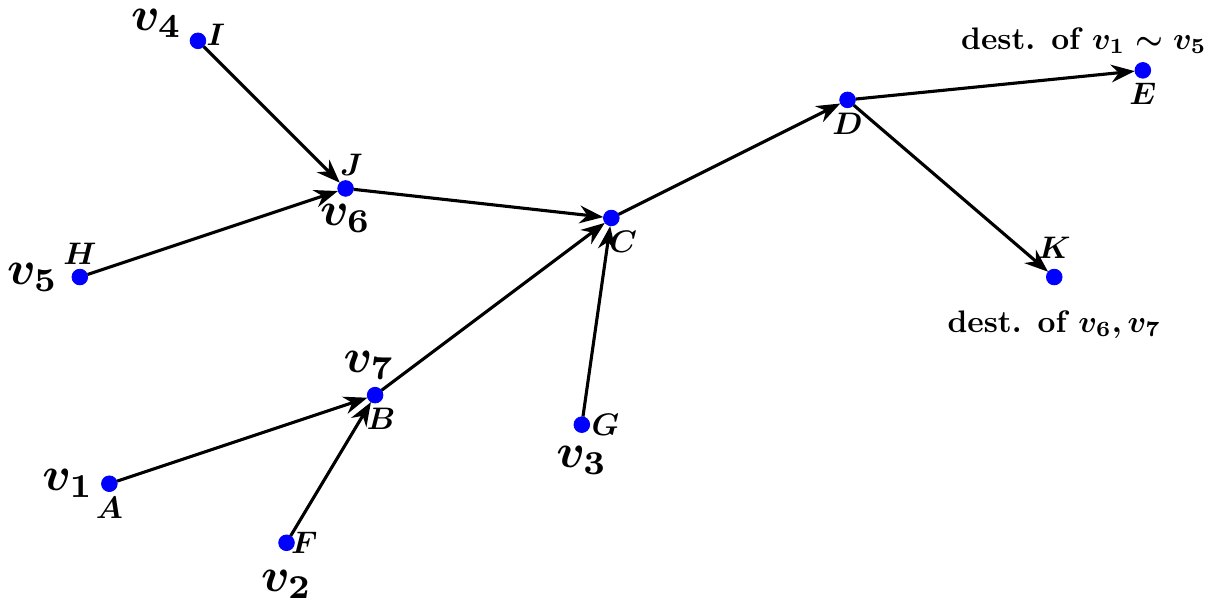}
    \caption{\scriptsize An illustrative figure of the polytree formed by the routes of 
    the vehicle set $\{v_1,v_2,v_3,v_4,v_5,v_6,v_7\}$. The departure nodes of the seven vehicles
    are $A$, $F$, $G$, $I$, $H$, $J$ and $B$, respectively. Vehicles $v_1\sim v_5$ have a shared
    destination node $E$, and the vehicles $v_6$ and $v_7$ have a shared destination node $K$. 
     Their routes are given as above and the direction of each arc is defined naturally by the travel
     direction of each vehicle along their routes. 
     }
    \label{fig:maximal-ideal-path}
\end{figure}

\noindent An example of a polytree formed by the union of vehicles' routes 
is given in Figure~\ref{fig:maximal-ideal-path}. 
\obsref{obs:plat-cond} gives a necessary and sufficient
condition on whether a subset of vehicles 
can form a platoon on their shared dir-path.

\begin{observation}\label{obs:plat-cond}
Let $P$ be a dir-path in $\cG$ shared by a vehicle set $\cV^{\prime}$. 
Let $s$ be the first node of $P$. 
By coordinating the departure time of vehicles in $\cV^\prime$, 
they can form a single platoon on $P$ (assuming no capacity constraint on a platoon) 
if and only if the intersection interval 
$\cap_{v\in\cV^\prime}[\ul{t}_{v,s},\ol{t}_{v,s}]$ is non-empty.
\end{observation}

\subsection{Transformation of time windows}
\label{sec:mst-rtw}
\begin{figure}
\centering
    \includegraphics[width=0.5\linewidth,trim=4cm 16.5cm 3.5cm 4.5cm]{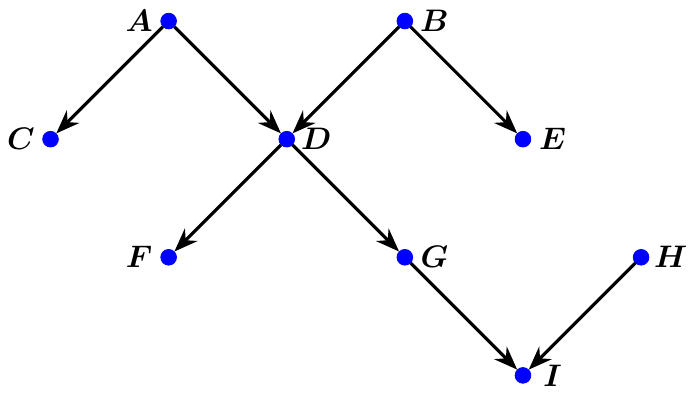}
    \caption{\scriptsize An illustrative figure of poly-paths in a polytree $\cG$. 
     The set $P^{\cG}_{C,E}=\{(A,C), (A,D), (B,D), (B,E)\}$ of arcs form a poly-path between $C$ and $E$,
     with $P^{\cG}_{C\to{E}}=\{(A,D),(B,E)\}$ and $P^{\cG}_{E\to{C}}=\{(B,D),(A,C)\}$.
     The set $P^{\cG}_{C,H}=\{(A,C),(A,D),(D,G),(G,I),(H,I)\}$ of arcs form a poly-path between $C$ and $H$
     with $P^{\cG}_{C\to{H}}=\{(A,D),(D,G),(G,I)\}$ and $P^{\cG}_{H\to{C}}=\{(H,I),(A,C)\}$. }
    \label{fig:poly-path}
\end{figure}

For the polytree $\cG$, 
we can select an arbitrary node $r$ as the \textit{reference node} of $\cG$. 
We can define a relative time window (RTW) for each vehicle $v\in\cV$
with respect to the reference node, which is a shift of the time window $[\ol{t}_{v,O_v},\ul{t}_{v,D_v}]$. 
With the help of this transformation, 
the platoonability of vehicles on their shared paths can be 
determined by checking if the intersection of their RTW's is non-empty.  
We first introduce a notion of poly-path which will be used in defining relative time difference. 
\begin{definition}[poly-path]
Let $A$ and $B$ be any two nodes of a polytree $\cG$ (with $\cE$ be the set of arcs).
The \textit{poly-path} between $A$ and $B$ is a subset $P^{\cG}_{A,B}=\{e_i\in\cE:i\in{I}\}$ 
($I$ being some set of indices) such that the set of undirected edges $\{e^{\textrm{undir}}_i\in\cE:i\in{I}\}$
form an undirected path between $A$ and $B$, where $e^{\textrm{undir}}_i$ represents the undirected
counterpart of $e_i$. Observe that every arc from $P^{\cG}_{A,B}$ (or equivalently $P^{\cG}_{B,A}$) 
is either in the direction from $A$ to $B$ ($A\to{B}$) or in the direction from $B$ to $A$ ($B\to{A}$). 
Denote by $P^{\cG}_{A{\to}B}$ the subset of arcs
from $P^{\cG}_{A,B}$ that is in the direction $A\to{B}$, and define $P^{\cG}_{B{\to}A}$ similarly. 
(An example of poly-paths in a polytree is given in \figref{fig:poly-path}.)
\end{definition}

\begin{definition}[relative time difference, relative departure time, relative time window]
\label{def:RTW}
Let $r$ be the reference node of the polytree $\cG$ (with $\cE$ being the set of arcs), 
and $s$ be any node of $\cG$.
(i) The \textit{relative time difference} of $s$ with respect to the reference node $r$ of $\cG$ is   
\begin{equation}
  \rT_{s,r}:=\sum_{(i,j)\in P^{\cG}_{s\to{r}}}T_{i,j}-\sum_{(i,j)\in P^{\cG}_{r\to{s}}} T_{i,j}.
\end{equation}  
(ii) Let $t$ be a feasible departure time of vehicle $v$ at its origin $O_v$. The \textit{relative departure time}
with respect to $r$ is defined as $t^{rel} = t + \rT_{O_v,r}$, where $\rT_{O_v,r}$ is the relative time difference 
of $O_v$ with respect to $r$.
(iii) For a $v\in\cV$, 
the \textit{relative time window} (RTW) of $v$ with respect to $r$
is defined as the interval: $\rtw_v:=[\underline{t}_{v,O_v}+\rT_{O_v,r},\; \overline{t}_{v,O_v}+\rT_{O_v,r}]$.
\end{definition}
In the above definition, the relative time difference $\rT_{s,r}$ can be viewed as a time shift
from node $s$ to node $r$ based on the poly-path linked them. 
The key motivation of introducing the notion of relative departure time and the relative time window
is to develop a relative time coordinate with respect to a reference node, such that vehicles
sharing a dir-path can form a platoon on it if and only if their relative departure times 
(with respect to the reference node) are equal. To provide an intuition,
consider the route graph $\cG$ given in \figref{fig:maximal-ideal-path} and
suppose the traveling time of each arc is 1. Focus on two
vehicles $\{v_1,v_3\}$. 
Consider the following departure times: $t_{v_1,A}=5$ and $t_{v_3,G}=6$.
The two vehicles enter the node $C$ simultaneously at time $7$, and hence
they can join a platoon on arc $CD$. Note that $t_{v_1,A}\neq{t}_{v_3,G}$ in this case.
However, if we let $A$ be the reference node and transform their departure times 
with respect to $A$. Their relative departure times are 
$t^\textrm{rel}_{v_1,A}=t_{v_1,A} + \rT_{A,A}=t_{v_1,A}=5$
and $t^\textrm{rel}_{v_3,G}=t_{v_3,G} + \rT_{G,A}=6 + (1-2)=5$,
and hence $t^\textrm{rel}_{v_1,A}=t^\textrm{rel}_{v_3,G}$.
Intuitively, it seems that two vehicles are platoonable on their shared arcs
if and only if their relative departure times are equal. If this property holds,
we can then transform the departure time windows of all vehicles to put them
into a relative-time based coordinate and discretize them in an efficient manner
such that if vehicles are assigned to a same time bucket, they are platoonable
on all of their shared arcs. This is the key idea of developing a formulation
without using big-M coefficients.

\begin{observation}\label{obs:single-path}
Let $\cV$ be a set of vehicles, and the route graph $\cG=\cup_{v\in\cV}\cR_v$ is a polytree.
For any $u,v\in\cV$, either $\cR_u\cap\cR_v=\emptyset$ or $\cR_u\cap\cR_v$
is a single directed path in $\cG$.
\end{observation}

\begin{example}
We given an example for the definition of relative time difference and relative time window
based on Figure~\ref{fig:maximal-ideal-path}. Suppose $A$ is chosen to be the reference node. 
The relative time difference of node $I$ with respect to node $A$ is given as
\bdm
\rT_{I,A} = T_{I,J}+T_{J,C}-T_{B,C}-T_{A,B},
\edm
where $P^{\cG}_{I\to{A}}=\{(I,J),\;(J,C)\}$ and $P^{\cG}_{A\to{I}}=\{(B,C),\;(A,B)\}$ in this case.
The RTW of $v_4$ is given by the interval $\rtw_{v_4}=[\ut_{v_4,I}+\rT_{I,A},\;\ot_{v_4,I}+\rT_{I,A}]$.
\end{example}

For the polytree $\cG$ with the reference node $r$,
we can calculate the RTW of each vehicle in $\cV$ following \defref{def:RTW}. 
\begin{proposition}\label{prop:pseudo-plat}
Suppose \assref{ass:polytree} hold. The following two properties hold: 
\emph{(a)} For any two vehicles $u,v\in\cV$ with $\cR_u\cap\cR_v\neq\emptyset$,
they can be platooned on $\cR_u\cap\cR_v$ if the departure times $t_{u,O_u}$, $t_{v,O_v}$
satisfy the following equation:
\begin{equation}
t_{u,O_u}+\rT_{O_u,r}=t_{v,O_v}+\rT_{O_v,r},
\end{equation}
i.e., $t^{rel}_{u,O_u}=t^{rel}_{v,O_v}$, where $r$ is the reference node of $\cG$.  
\end{proposition}
\proof{Proof.}
By \obsref{obs:single-path}, $\cR_u\cap\cR_v$ is a directed path.
Let $s$ be the first node on $\cR_u\cap\cR_v$. 
A necessary and sufficient condition for $u$ and $v$
to form a platoon on the shared path is that the arrival times of $u$ and $v$
at the node $s$ must be equal.
Let $\cR_u(O_u,s)$ be the 
directed path from $O_u$ to $s$, and let $\cR_v(O_v,s)$ be the directed path from
$O_v$ to $s$. We must have $\cR_u(O_u,s)\cap\cR_v(O_v,s)=\emptyset$.
The simultaneous arrival condition for $u$ and $v$ at node $s$ is written as:
\begin{equation}\label{eqn:t+a=t+a}
  t_{u,O_u}+\sum_{(i,j)\in\cR_u(O_u,s)}T_{i,j}=t_{v,O_v}+\sum_{(i,j)\in\cR_v(O_v,s)}T_{i,j}.
\end{equation}
By definition, the relative time differences of $O_u$ and $O_v$ with respective to $r$ are
\begin{equation}\label{eqn:T=T}
\begin{aligned}
  &\rT_{O_u,r}=\sum_{(i,j)\in\cR_u(O_u,s)}T_{i,j}+\sum_{(i,j)\in P^\cG_{s\to r}}T_{i,j}-\sum_{(i,j)\in P^\cG_{r\to s}}T_{i,j}, \\
  &\rT_{O_v,r}=\sum_{(i,j)\in\cR_v(O_v,s)}T_{i,j}+\sum_{(i,j)\in P^\cG_{s\to r}}T_{i,j}-\sum_{(i,j)\in P^\cG_{r\to s}}T_{i,j},
\end{aligned}
\end{equation} 
where the notations $P^\cG_{s\to r}$ and $P^\cG_{r\to s}$ are given in \defref{def:RTW}. 
Using \eqref{eqn:T=T}, we find that the equation \eqref{eqn:t+a=t+a} is equivalent to 
$t_{u,O_u}+\rT_{O_u,r}=t_{v,O_v}+\rT_{O_v,r}$, concluding the proof.   
\halmos
\endproof
     
In the following example, we illustrate how \propref{prop:pseudo-plat}
can lead to an efficient time discretization scheme and a MILP formulation of the tree-FRCVP
without using big-M coefficients. 
\begin{example}\label{exp:interval-graph}
For the vehicle routes given in Figure~\ref{fig:maximal-ideal-path},
suppose the travel time of each arc is 1, and
we consider the following departure time and arrival time restrictions:
the $(T^O_v,T^D_v)$ for $v_1\sim v_7$ is $(4,12)$, $(3,11)$, $(4,13)$,
$(4,11)$, $(9,19)$, $(11,15)$ and $(13,19)$, respectively.
Let $A$ be the reference node of the polytree. Using Definition~\ref{def:time-wind},
we can obtain the RTW's of vehicles $v_1\sim v_7$ as $[4,8]$,
$[3,7]$, $[3,9]$, $[4,7]$, $[9,15]$, $[10,11]$ and $[12,15]$, respectively.
The tree-FRCVP problem can be viewed as a coordination problem
of their relative departure times constrained by the RTW's to let as many
vehicles having the same relative departure times. 
Now consider vehicles $\{v_1,v_2,v_3\}$. 
To ensure they are platoonable on the shared arcs, we need to select
a same relative departure time $t^{rel}$ for them according to \propref{prop:pseudo-plat}.  
Therefore, $t^{rel}$ can only be selected from the intersection of their RTW's,
i.e., $t^{rel}\in[4,8]\cap[3,7]\cap[3,9]=[4,7]$. A key observation is that the platooning configurations
among the three vehicles on their shared arcs are independent of the specific value
of $t^{rel}$ as long as $t^{rel}\in[4,7]$. This implies that we can transform the original 
continuous-time coordination problem into a discrete-time coordination problem
by generating all possible time intervals (or time buckets) yielded by all the RTW's.
Then the decision-making is about assigning each vehicle to a time bucket as a 
representation of its relative departure time.
\end{example}

\subsection{Efficient time discretization and a vehicle-to-time-bucket assignment formulation}
\label{sec:time-discret}
\begin{figure}
\centering
    \includegraphics[width=0.47\linewidth,trim=2.5cm 18cm 2.5cm 4cm]{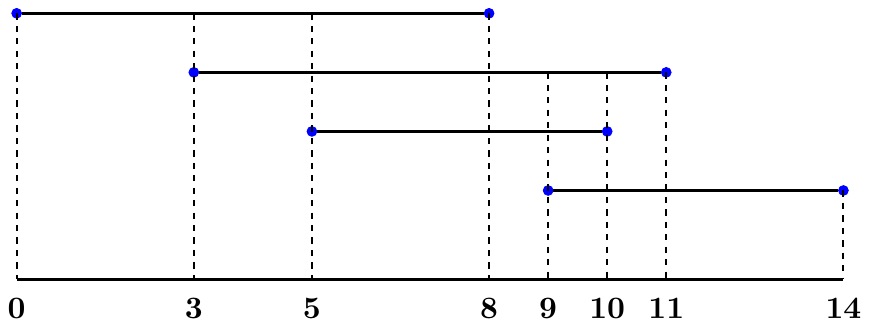}
    \caption{The time buckets generated by \algoref{alg:ATD} when applied to the four-vehicle 
    instance given in Example~\ref{exp:ATD}.}
    \label{fig:ATD}
\end{figure}

In this section, we formally develop an efficient time-discretization (ETD) method that leads to an improved
formulation for the tree-FRCVP problem, based on the result of \propref{prop:pseudo-plat}. 
The essence of ETD is to generate time buckets based on
the intersected and non-intersected segments of all RTW's,
by leveraging the result of \propref{prop:pseudo-plat}.
The ETD method is presented in Algorithm~\ref{alg:ATD},
and it is illustrated in Example~\ref{exp:ATD} and Figure~\ref{fig:ATD}(b). 

\begin{definition}[sorted time buckets]
A set $\cS=\{[p_k,q_k]\}^{|\cS|}_{k=1}$ is a \emph{collection of sorted time buckets}
if $\cS$ is finite and $q_k\le p_{k+1}$ for all $k\in\{1,\ldots,|\cS|-1\}$. 
\end{definition}

\begin{algorithm}
{\scriptsize 
\caption{\footnotesize Efficient time discretization.}\label{alg:ATD}
\begin{algorithmic}[1]
\State{\bf Input}: a tree-FRCVP instance with the vehicle set $\cV=\{1,\ldots,N\}$.
Denote $[\tau^k_1,\tau^k_2]$ as the RTW for each $k\in\cV$.  
\State{\bf Output}: A set $\cS=\{[p_k,q_k]\}^{|\cS|}_{k=1}$ of sorted time buckets.
\State{Set $\cS=\emptyset$, $\mL\gets\textrm{ empty list}$}.
\For{$k\in\cV$}
	\State{$\mL\gets\mL.append\{\tau^k_1,\tau^k_2\}$ (duplication of elements is allowed)}
\EndFor
\State{Sort the list $\mL$ (duplicated elements are consecutive in the sorted list).}
\For{$k=1,\dots,|\mL|-1$}
	\State{Get the $k$-th and $(k+1)$-th elements of $\mL$, and denote them as $a_k$ and $a_{k+1}$, respectively.}
	\If{$[a_k,a_{k+1}]$ is a subset of $[\tau^m_1,\tau^m_2]$ for some $m\in\cV$}
		\State{Append the time bucket $[a_k,a_{k+1}]$ to the end of $\cS$.}
	\EndIf
\EndFor
\State{\Return{$\cS$}.}
\end{algorithmic}
}
\end{algorithm}
\begin{remark}
The reason of allowing duplicated elements in the list $\cL$ in \algoref{alg:ATD} is demonstrated
by the following corner case: Consider two vehicles with RTW's being $[1,3]$ and $[3,5]$, respectively.
In the continuous-time setting, both vehicles can select 3 as their common relative time and 
determine their departure times accordingly. This decision ensures that they are eligible to form a
platoon in the shared dir-path. Allowing duplicated elements in $\cL$ can lead to generating time buckets
$[1,3]$, $[3,3]$ and $[3,5]$. Then $[3,3]$ is a time bucket feasible to both vehicles, and they can be 
assigned to $[3,3]$ to form a platoon in shared arcs.
\end{remark}

\begin{example}\label{exp:ATD}
We consider four relative time windows (RTWs) given as $[0,8]$, $[3,11]$, $[5,10]$ and $[9,14]$. 
If Algorithm~\ref{alg:ATD} is applied to the above
RTWs, it will generate 7 sorted time buckets: $[0,3]$, $[3,5]$, $[5,8]$,
$[8,9]$, $[9,10]$, $[10,11]$ and $[11,14]$ (\figref{fig:ATD}(b)).
\end{example}

\begin{restatable}{observation}{propIII}\label{prop:determ-bound-|S|} 
Let $\cS$ be the set of sorted time buckets returned by \algoref{alg:ATD}. 
Then the relation $1\le|\cS|\le 2|\cV|-1$ holds.
\end{restatable}

\begin{table}
  \footnotesize
  \caption{Variables in the vehicle-to-time-bucket assignment formulation \eqref{opt:VA} for the tree-FRCVP problem. \label{tab:var-TWOF}}
  \begin{tabularx}{\textwidth}{lX}
    \toprule
    Set & Definition \\
    \midrule
		$\cS$ & the set of time buckets, i.e., the output of Algorithm~\ref{alg:ATD}; \\
		$\cS_u$ & the subset of time buckets that are feasible to vehicle $u$, i.e., $\cup_{t\in\cS_u}t=\text{RTW}_u$, $u\in\cV$; \\
		$\cV_{e,t}$ & the subset of vehicles that share arc $e\in\cR$ and are feasible at time bucket $t\in\cS$; \\
    \midrule\\[-.35cm]
    Parameter & Definition \\
    \midrule
    		$\lambda\ge2$ & the capacity of a platoon, i.e., the maximal number of vehicles allowed in a platoon \\
    \midrule\\[-.35cm]
    Variable & Definition \\
    \midrule 
		$x_{u,t}$ & 1 if vehicle $u$ is assigned to time bucket $t\in\cS$, 0 otherwise;  \\
		$z_{e,t}$ & the number of full-size platoons (i.e., a platoon of size $\lambda$) taking $e\in\cR$ at time bucket $t\in\cS$; \\
		$y_{e,t}$ & 1 if there exists a unsaturated platoon (i.e., a platoon of size at least one but less than $\lambda$)
				  taking $e\in\cR$ at time bucket $t\in\cS$, 0 otherwise; \\
		$y^\prime_{e,t}$ & 1 if the unsaturated platoon taking $e\in\cR$ at time bucket $t\in\cS$ is of the size at least 2,
				 0 otherwise; \\
		$q_{e,t}$ & the size of the unsaturated platoon taking $e\in\cR$ at time bucket $t\in\cS$; \\
		$w_{e,t}$ & defined as $\sum_{u\in\cV_{e,t}}x_{u,t}-z_{e,t}-y_{e,t}$, 
			which can be interpreted as the total number of trailing vehicles\footnotemark[2] of platoons taking $e\in\cR$ at $t\in\cS$.  \\  
    \bottomrule
  \end{tabularx}
\end{table}
\footnotetext[2]{Given $\lambda=5$, if there are 11 vehicles taking $e\in\cR$ at $t\in\cS$, 
they form 3 platoons: 2 saturated platoons and 1 single-vehicle platoon (unsaturated), and then $w_{e,t}=11-3=8$.
If there are 12 vehicles taking $e\in\cR$ at $t\in\cS$, they form 3 platoons: 2 saturated platoons and 1 unsaturated platoon of size 2,
and then $w_{e,t}=12-3=9$.}

The key observation pointed by the efficient time discretization
is that the tree-FRCVP problem can be formulated equivalently as an advanced assignment
problem. Every vehicle is allowed to be assigned 
to a set of feasible consecutive time buckets, and vehicles that are
assigned to a same time bucket are platoonable on their shared arcs. 
The assignment formulation is given as \eqref{opt:VA} 
with sets and variables defined in Table~\ref{tab:var-TWOF}.
\begin{subequations}
	\makeatletter
	\def\@currentlabel{VA}
	\makeatother
	\label{opt:VA}
	\renewcommand{\theequation}{VA.\arabic{equation}}
	\addtolength{\jot}{0.1mm}
		\begin{alignat}{3}
			\max  &\quad \sum_{e\in\cR}\sum_{t\in\cS}(\sigma^lC_ez_{e,t}
			+\sigma^lC_ey^\prime_{e,t}+\sigma^fC_ew_{e,t})    &  \label{opt:va1} \\
			\text{s.t. } & \sum_{t\in\cS_u}x_{u,t}=1 &\; \forall u\in\cV, \label{opt:va2} \\
			& \lambda z_{e,t}-\sum_{u\in\cV_{e,t}}x_{u,t} \le 0  &\; \forall  e\in\cR,\; t\in\cS,  \label{opt:va3} \\
			& \lambda (z_{e,t}+1)-\sum_{u\in\cV_{e,t}}x_{u,t} \ge 1 &\; \forall  e\in\cR,\; t\in\cS,  \label{opt:va4} \\
			& w_{e,t}-\sum_{u\in\cV_{e,t}}x_{u,t}+z_{e,t}+y_{e,t}= 0  &\; \forall e\in\cR,\; t\in\cS,  \label{opt:va5} \\
			& q_{e,t} = \sum_{u\in\cV_{e,t}}x_{u,t} - \lambda z_{e,t} &\; \forall  e\in\cR,\; t\in\cS,  \label{opt:va6} \\
			& y_{e,t} \le q_{e,t},\;\; \lambda y_{e,t} \ge q_{e,t} &\; \forall  e\in\cR,\; t\in\cS,  \label{opt:va7} \\
			& q_{e,t}\ge 2y^\prime_{e,t},\;y^\prime_{e,t}\le y_{e,t}  &\; \forall  e\in\cR,\; t\in\cS,  \label{opt:va8} \\ 	 
			& x_{u,t},\; y_{e,t},\; y^\prime_{e,t}\in\{0,1\},\;\; z_{e,t},\; q_{e,t},\; w_{e,t}\in\bZ_+ &\; \forall e\in\cR, \; t\in\cS,\; u\in\cV_{e,t}. \label{opt:va10}
		\end{alignat}
\end{subequations}
Vehicles in a platoon can be classified into three types based on the fuel saving: 
a \textit{lead vehicle}, a \textit{trailing vehicle} and a \textit{single vehicle}. In a platoon of size at least two, 
the lead vehicle is the one at the head of the platoon and trailing vehicles are the other vehicles in the platoon.
A single vehicle is the only one in a single-vehicle platoon which cannot achieve any fuel saving.
The objective is to maximize the total fuel saving of all vehicles.
The total fuel saving is a sum of the fuel saving achieved at every $(e,t)$
combination. For each $(e,t)$,
the fuel saving is determined by the number of lead vehicles $z_{e,t}+y^\prime_{e,t}$
and the number of trailing vehicles $w_{e,t}$. 
The constraint \eqref{opt:va2} ensures that each vehicle must be assigned 
to exact one time bucket that is feasible for it. The constraints \eqref{opt:va3} and \eqref{opt:va4}
ensure that $z_{e,t}$ is the quotient of the total number of vehicles at $(e,t)$ divided by
the platoon size limit $\lambda$. The constraint \eqref{opt:va5} rewrites the total number of
vehicles in terms of $w_{e,t}$, $z_{e,t}$ and $y_{e,t}$ for each $(e,t)$.
The constraint \eqref{opt:va6} ensures that $q_{e,t}$ is the number of vehicles
involved in the unsaturated platoon, and its value is zero if all platoons are saturated. 
The constraint \eqref{opt:va7} ensures that $y_{e,t}$ is equal to one if and only if $q_{e,t}\ge 1$.
The constraints in \eqref{opt:va8} are only necessary to define $y^\prime_{e,t}$. However,
when combining with the sense of maximization, $y^\prime_{e,t}$
will take the value 1 instead of 0 as long as the two constraints are respected. 
\begin{theorem}
An optimal solution of \eqref{opt:VA} can be transformed into an optimal solution 
of the corresponding tree-FRCVP problem. 
\end{theorem}
\proof{Proof.}
Let $(x^*,y^*,y^{\prime*},z^*,q^*,w^*)$ be an optimal solution of \eqref{opt:VA}. 
Let $\cV^*_t=\{v\in\cV: x^*_{v,t}=1\}$ and $\cV^*_{e,t}=\{v\in\cV: x^*_{v,t}=1,\; e\in\cR_v\}$.
For each time bucket $b\in\cS$, select an arbitrary value $\tau_b$ from the interval represented by $b$
and set the departure time $t_{v,O_v}=\tau_b-\wt{T}_{O_v,r}$ for every $v\in\cV^*_b$. 
The set of departure times $dpt:=\{t_{v,O_v}: v\in\cV\}$ induces a feasible solution of the
corresponding tree-FRCVP.
The platooning configuration induced by $dpt$ can be determined based on the simultaneousness
of entering time of related vehicles at their shared arcs. 
By \propref{prop:pseudo-plat}, for every $b\in\cS$ and $e\in\cE$, all vehicles from the subset $\cV^*_{e,b}$
enter the arc $e$ simultaneously according to $dpt$. Let vehicles in $\cV^*_{e,b}$ form saturated platoons
as many as the count $|\cV^*_{e,b}|$ permits and leave at most one unsaturated platoon.
Due to constraints \eqref{opt:va3}$\sim$\eqref{opt:va6},
they can form $z^*_{e,b}$ platoons of size $\lambda$ on arc $e$ with $q^*_{e,b}$ ($q^*_{e,b}<\lambda$) 
left over vehicles forming a platoon. This shows that the $dpt$ as a solution of the corresponding tree-FRCVP instance
has the same objective value as the \eqref{opt:VA}, i.e., $\obj(dpt)\ge\textproc{OPT}\eqref{opt:VA}$. 

It remains to show that $\obj(dpt^\prime)\le\textproc{OPT}\eqref{opt:VA}$ for any solution $dpt^\prime$
of the tree-FRCVP instance. Suppose $dpt^\prime=\{t^\prime_{v,O_v}: v\in\cV\}$ is any set of departure times
and $\tau^\prime=\{\tau^\prime_v: v\in\cV\}$ is the set of relative departure times derived from 
$\{t^\prime_{v,O_v}: v\in\cV\}$ satisfying $\tau^\prime_v=t^\prime_{v,O_v}+\wt{T}_{O_v,r}$.
Based on $\{\tau^\prime_v: v\in\cV\}$, we can construct a feasible assignment $x^\prime$ of \eqref{opt:VA}
by assigning $v$ to the time bucket $b$ that contains $\tau^\prime_v$. Note that if $\tau^\prime_v$
is at the boundary of two consecutive time buckets, assign it to the earlier bucket.  
In this construction, $\tau^\prime_u=\tau^\prime_v$ implies that $u$ and $v$ will be assigned into
a same time bucket, but the reverse statement does not hold. 
This shows that the assignment $x^\prime$ can lead to a platooning configuration 
as good as the one yielded by $dpt^\prime$. This shows that 
$\obj(dpt^\prime)\le\textproc{OPT}\eqref{opt:VA}$, which concludes the proof.
\halmos
\endproof


\tabref{tab:form-compare} summarizes the comparison between 
the continuous-time formulation and vehicle-to-time-bucket assignment
in terms of the complexity and formulation power.
\begin{definition}
Consider the set $\cV$ of vehicles involved in a polytree.
A time bucket is \emph{feasible} to vehicle $v$ if it is contained in $\rtw_v$.
\end{definition}
The quantity $|\cS_u|$ is the number of feasible time buckets of the vehicle $u$,
and it is a crucial factor that determines the complexity of the formulation \eqref{opt:VA}. 
The following proposition gives a result regarding the expectation and 
a tail bound of $|\cS_u|$ when the RTW of each vehicle is determined by i.i.d. random parameters.
The proof is given in \ref{sec:suppl-1}.
\begin{restatable}{proposition}{thmI}\label{thm:E[|S|]}
Suppose \assref{ass:polytree} hold, and suppose 
for any vehicle $v\in\cV$ the $\rtw_v$ has the form $[T_v,T_v+G_v]$,
where $\{T_v:v\in\cV\}$ and $\{G_v:v\in\cV\}$ are i.i.d. random variables 
following continuous distribution $F_T$ and $F_G$, respectively, and the 
two sets of random variables are also independent. Then for any $u\in\cV$,
\begin{subequations}
\begin{equation}\label{eqn:E[Su]}
\E[|\cS_u|] = 1 + (|\cV|-1)\mu, 
\end{equation}
\begin{equation}\label{eqn:P(Su)}
P\left(\left|\frac{|\cS_u|-1}{|\cV|-1}-\mu\right|\ge t\right) \le 2\exp\left(-\frac{(|\cV|-1)t^2}{2}\right),
\end{equation}
\end{subequations}
where 
\beq\label{eqn:mu}
\ba
\mu&=\frac{1}{2}P(\min\{G_1,G_2\}>|T_1-T_2| ) + \frac{1}{2}P(\max\{G_1,G_2\}>|T_1-T_2|>|G_1-G_2|) \\
&\quad+ \frac{1}{2}P(|G_1-G_2|>|T_1-T_2|),
\ea
\eeq
$T_i$ and $G_i$ ($i=1,2$) are independent random variables satisfying $T_1,T_2\sim F_T$ 
and $G_1,G_2\sim F_G$.
\end{restatable}
\begin{restatable}{corollary}{corUniform}\label{cor:uniform}
If $T\sim\textrm{Uniform}[0,H]$ and $G\sim\textrm{Uniform}[0,h]$ for $0<h<H$, 
then $$\mu=\frac{h}{H}-\frac{h^2}{3H^2}.$$
\end{restatable}

\begin{table}
  \scriptsize 
  \caption{\scriptsize Comparison of two formulations for the tree-FRCVP problem.}  \label{tab:form-compare}
  \begin{tabularx}{\textwidth}{lllll}
    \toprule
    Formulation & Variables & Constraints & big-M coefficients & Network topology  \\
    \midrule
    \eqref{opt:cont_time} continuous-time not allowing stop-and-wait & $\mathcal{O}(|\cV|^2|\cE|)$ & $\mathcal{O}(|\cV|^2|\cE|)$ & Yes & any \\
    \eqref{opt:VA} discrete-time not allowing stop-and-wait & $\mathcal{O}{(|\cV||\cS|+|\cE||\cS|)}$ & $\mathcal{O}{(|\cE||\cS|)}$ & No & tree \\
    \bottomrule
  \end{tabularx}
\end{table}

\section{Discrete-time formulation for tree-FRCVP instances allowing waiting at intermediate nodes}
\label{sec:allow-waiting}
In this section, we analyze for the tree-FRCVP instances in which the stop-and-wait option 
at intermediate nodes are permitted, i.e., \assref{ass:no-pause} does not hold.
This means each vehicle $v$ has the flexibility to stop at any intermediate nodes ($\cN_v\setminus\{O_v,D_v\}$) 
along $\cR_v$ and wait for a moment to increase the opportunities of engaging in different platoons to save more energy.
There is no limit on the number of stopping times and the amount of waiting time as long as
the origin and destination time restrictions are respected.
We refer this problem as  \textbf{tree-FRCVP-pause} problem.

Each intermediate stop may incur a fixed amount of engine-setup fuel cost $\alpha_{\textrm{eng}}$
when restarting the vehicle. 
For the tree-FRCVP-pause problem, there are 
additional decisions to make: For each $v\in\cV$ and $s\in\cN_v\setminus\{O_v,D_v\}$, decide whether
$v$ should wait at the node $s$ and how much time it should wait.  
We will show that there is a (big-M free) discrete-time MILP formulation 
for the tree-FRCVP-pause problem based on extending \eqref{opt:VA}.
In particular, we consider the following modification of decision variables:
\begin{enumerate}
	\item Generate feasible set of time buckets $\cS_u$ for each vehicle $u\in\cV$ using \algoref{alg:ATD} in the same way as in \secref{sec:time-discret}.
	\item Create vehicle-arc pairs $(u,e)$ for every $u\in\cV$ and every $e\in\cR_u$.	  
	\item Decisions: For every vehicle-arc pair $(u,e)$, assign it to exactly one of the feasible time buckets from $\cS_u$.
	\item Additional constraints: If an arc $e$ is in front an arc $e^\prime$ along $\cR_u$, 
		the time bucket $t$ that $(u,e)$ has been assigned cannot be after the time bucket
		$t^\prime$ that $(u,e^\prime)$ has been assigned. (But $(u,e)$ and $(u,e^\prime)$ can be assigned to a same time bucket.)
\end{enumerate}
For simplicity we consider the case $\lambda=\infty$ (no limitation on the size of a platoon). 
The discrete time MILP formulation for the tree-FRCVP-pause problem is given in {\DT} 
with new parameters and decision variables defined in \tabref{tab:var-DT}.
The formulation {\DT} manages to model the following logic: 
For any two consecutive arcs $e=(n_1,n_2)$
and $e^\prime=(n_2,n_3)$ along $\cR_u$, if $(u,e)$ is assigned to a time bucket $t$
but $(u,e^\prime)$ is assigned to a different bucket $t^\prime$ (with $t<t^\prime$ given feasibility),
if and only if $u$ has stopped at node $n_2$ and wait for a positive amount of time.

\begin{table}
  \footnotesize
  \caption{ Variables in the discrete-time formulation {\DT} for the tree-FRCVP-pause problem. } \label{tab:var-DT}
  \begin{tabularx}{\textwidth}{lX}
  \toprule
     Parameter & Definition \\
     \midrule
    		$\alpha_{\textrm{eng}}$ & the energy cost for a vehicle to restart the engine after waiting at an intermediate node; \\
    \midrule\\[-.35cm]
    Variable & Definition \\
    \midrule 
		$x_{u,e,t}$ & 1 if vehicle $u$ is assigned to time bucket $t\in\cS$ for traversing the arc $e$, 0 otherwise;  \\
		$p_{u,t}$ & equals $\max_{e\in\cR_u}x_{u,e,t}$; \\ 
		$y_{e,t}$ & 1 if there is at least one vehicle assigned to the time bucket $t$ for traversing $e\in\cR$, 0 otherwise; \\
		$y^\prime_{e,t}$ & 1 if there are at least two vehicles assigned to the time bucket $t$ for traversing $e\in\cR$, 0 otherwise; \\
		$w_{e,t}$ & $N_{e,t}-1$ if $N_{e,t}\ge1$, 0 otherwise. $N_{e,t}$ is the number of vehicles assigned to $t$ 
				   for traversing $e\in\cR$.  \\  
    \bottomrule
  \end{tabularx}
\end{table}

\begin{subequations}
	\makeatletter
	\def\@currentlabel{DT}
	\makeatother
	\label{opt:DT}
	\renewcommand{\theequation}{DT.\arabic{equation}}
	\addtolength{\jot}{0.1mm}
		\begin{alignat}{3}
		\max  &\quad \sum_{e\in\cR}\sum_{t\in\cS}(\sigma^lC_ey^\prime_{e,t}+\sigma^fC_ew_{e,t})-\alpha_{\textrm{eng}}\sum_{u\in\cV}\Big(\sum_{t\in\cS_u}p_{u,t}-1 \Big)   &  \label{opt:dt-obj}  \\
			\text{ s.t. } & \sum_{t\in\cS_u}x_{u,e,t}=1 & \forall u\in\cV,\; e\in\cR_u, \label{opt:dt-c1} \\
			& x_{u,e^\prime,s}\le 1-x_{u,e,t} & \forall e\prec{e}^\prime\in\cR_u,\;s\prec{t}\in\cS_u, \label{opt:dt-c2} \\
			& p_{u,t} \ge x_{u,e,t} & \forall u\in\cV,\; e\in\cR_u,\;t\in\cS_u, \label{opt:dt-c3} \\
			& \sum_{u\in\cV_e}x_{u,e,t}\ge2y^\prime_{e,t} & \forall e\in\cR,\;t\in\cS, \label{opt:dt-c4} \\
			& y_{e,t}\ge x_{u,e,t} & \forall e\in\cR,\; u\in\cV_e,\;t\in\cS_u, \label{opt:dt-c5} \\
			& w_{e,t} = \sum_{u\in\cV_e}x_{u,e,t}-y_{e,t} & \forall e\in\cR,\;t\in\cS, \label{opt:dt-c6} \\
			& x_{u,e,t}\in\{0,1\} & \forall u\in\cV,\; e\in\cR_u,\; t\in\cS_u, \\
			& p_{u,t}\in\{0,1\} & \forall u\in\cV,\; t\in\cS_u, \\
			& y_{e,t},\;y^\prime_{e,t}\in\{0,1\},\;w_{e,t}\ge0 & \forall e\in\cR,\; t\in\cS,
		\end{alignat}
\end{subequations}
where the first term in the objective models the fuel saving achieved at every arc and by every platoon,
the second term is the extra fuel cost of restarting a vehicle if it stops at intermediate nodes. 
Note that $\sum_{t\in\cS_u}p_{u,t}$ is total number of time buckets from $\cS_u$ containing
at least one $(u,e)$ pair $e\in\cR_u$. This quantity equals one if and only if all $(u,e)$ pairs
are assigned to a single time bucket which is equivalent to that $u$ has not stopped on any 
intermediate node, otherwise $\sum_{t\in\cS_u}p_{u,t}-1$ is the number of intermediate stops
of $u$. \eqref{opt:dt-c1} ensures that each $(u,e)$ pair is assigned to exact one
time bucket from $\cS_u$. \eqref{opt:dt-c2} ensures the restriction that if $e^\prime$ is after $e$
along $\cR_u$, then $(u,e^\prime)$ cannot be assigned to any time bucket earlier than the time bucket
$(u,e)$ has been assigned to, where $\prec$ is the preceding relation describing the order of 
two arcs along a same route or two time buckets. \eqref{opt:dt-c3}$\sim$\eqref{opt:dt-c6} 
(together with the sense of maximization) translate
the definitions of $p_{u,t}$, $y^\prime_{e,t}$, $y_{e,t}$ and $w_{e,t}$, respectively.

Observe that a feasible solution of the tree-FRCVP-pause can be sufficiently
characterized using the the departure time $\hat{\tau}_{v,e}$ from the beginning node\footnotemark[3] of the arc $e$ for
every $e\in\cR_v$ and $v\in\cV$. If $e$ is the first arc on $\cR_v$ for a vehicle $u$, $\hat{\tau}_{v,e}$ gives the
departure time of $v$ from the origin node $O_v$. Otherwise, if $e=(n_1,n_2)$ is an arc from $\cR_v$ with $n_1\neq{O}_v$,
$\hat{\tau}_{v,e}$ gives the departure time of $v$ from the intermediate node $n_1$. The 
waiting time at node $n_1$ can be inferred based on the departure time $\hat{\tau}_{v,e^\prime}$ where
$e^\prime=(n_0,n_1)$ is the arc right before $e$ on $\cR_u$: If $\hat{\tau}_{v,e}>\hat{\tau}_{v,e^\prime}+T_{n_0,n_1}$,
then $v$ stops at $n_1$ and wait for $\hat{\tau}_{v,e}-\hat{\tau}_{v,e^\prime}-T_{n_0,n_1}$ amount of time.
The platoon formation on an arc $e$ can be inferred straightforwardly based on the 
departure times of vehicles (sharing $e$) from the beginning node of $e$:
All vehicles sharing an arc $e=(n_1,n_2)$ and departing from $n_1$ at the same (instant) time should form a platoon on $e$
subject to the platoon capacity.
Given a feasible solution $(\hat{x},\hat{p},\hat{y},\hat{y}^\prime)$ of {\DT},
it can be transformed into a feasible solution of the corresponding tree-FRCVP-pause instance using
\algoref{alg:sol-DT-transform} which is proved in \propref{prop:sol-transform}. 
\footnotetext[3]{For an arc $e=(n_1,n_2)$, $n_1$ is the beginning node of $e$.}

\begin{algorithm}
\footnotesize 
\caption{\footnotesize Transformation from a feasible solution of {\DT} to a feasible solution of tree-FRCVP-pause.}\label{alg:sol-DT-transform}
\begin{algorithmic}[1]
\State{\bf Input}: A feasible solution $(\hat{x},\hat{p},\hat{y},\hat{y}^\prime,\hat{w})$ of {\DT}.
\State{\bf Output}: A feasible solution $\{\hat{\tau}_{u,e}:u\in\cV,\;e\in\cR_u\}$, 
 where $\hat{\tau}_{u,e}$ is the departure time of $u$ from the starting node of $e$.
\For{$t:=[\ul{t},\ol{t}]\in\cS$}
	\State{Consider the set $\cF_t:=\{(u,e):\hat{x}_{u,e,t}=1\}$.}
	\For{$(u,e)\in\cF_t$}
		\State{Suppose $e=(n_1,n_2)$.}
		\State{Set $\hat{\tau}_{u,e}=\ul{t}-\wt{T}_{n_1,r}$, where $r$ is the reference node.}
	\EndFor
\EndFor
\State{\bf Return} $\{\hat{\tau}_{u,e}:u\in\cV,\;e\in\cR_u\}$.
\end{algorithmic}
\end{algorithm} 

\begin{proposition}\label{prop:sol-transform}
For a feasible solution $(\hat{x},\hat{p},\hat{y},\hat{y}^\prime,\hat{w})$ of \eqref{opt:DT},
the solution $\{\hat{\tau}_{u,e}:u\in\cV,\;e\in\cR_u\}$ returned by \algoref{alg:sol-DT-transform}
is feasible to the tree-FRCVP-pause instance, 
where $\hat{\tau}_{u,e}$ is the departure time of $u$ from the starting node of $e$.
\end{proposition}
\proof{Proof.}
For any $u\in\cV$ and $e=(n_1,n_2)\in\cR_u$, suppose $t=[\ul{t},\ol{t}]$
is the time bucket that $(u,e)$ has been assigned to by the feasible solution of \eqref{opt:DT}.\
The algorithm sets $\hat{\tau}_{u,e}=\ul{t}-\wt{T}_{n_1,r}$. 
Since $\ul{t}\in[\ul{t}_{u,O_u}+\wt{T}_{O_u,r},\ol{t}_{u,O_u}+\wt{T}_{O_u,r}]$,
one has 
\bdm
\ba
&\hat{\tau}_{u,e_0}\ge\ul{t}_{u,O_u}+\wt{T}_{O_u,r}-\wt{T}_{n_1,r}=\ul{t}_{u,O_u}+\sum_{e\in\cR_u(O_u,n_1)}T_e, \\
&\hat{\tau}_{u,e_0}\le\ol{t}_{u,O_u}+\wt{T}_{O_u,r}-\wt{T}_{n_1,r}=\ol{t}_{u,O_u}+\sum_{e\in\cR_u(O_u,n_1)}T_e,
\ea
\edm
where $\cR_u(O_u,n_1)$ is the set of arcs from $O_u$ to $n_1$ on $\cR_u$.
This shows that $\hat{\tau}_{u,e}\in[T^O_u,T^D_u]$ by the definition
 the departure time window of $u$ (\defref{def:time-wind}).

It remains to verify that for any $u\in\cV$ and two arcs $e_1\prec{e_2}\in\cR_u$
with $e_1=(n_1,n_2)$ and $e_2=(n_3,n_4)$,
we must have $\hat{\tau}_{u,e_2}\ge\hat{\tau}_{u,e_1}+\sum_{e\in\cR_u(n_1,n_3)}T_{e}$,
where $\cR_u(n_1,n_3)$ is the set of arcs from $n_1$ to $n_3$ on $\cR_u$.
Let $t_1=[\ul{t}_1,\ol{t}_1]$ and $t_2=[\ul{t}_2,\ol{t}_2]$ be the two time buckets satisfying
$\hat{x}_{u,e_1,t_1}=1$ and $\hat{x}_{u,e_2,t_2}=1$. The feasibility implies that $t_1\preceq{t_2}$
(either $t_1=t_2$ or $t_1\prec{t_2}$). By the algorithm, we have $\hat{\tau}_{u,e_1}=\ul{t}_1-\wt{T}_{n_1,r}$
and $\hat{\tau}_{u,e_2}=\ul{t}_2-\wt{T}_{n_3,r}$. Therefore, the following inequalities hold:
\bdm
\ba
\hat{\tau}_{u,e_2}-\hat{\tau}_{u,e_1}=(\ul{t}_2-\wt{T}_{n_3,r}) - (\ul{t}_1-\wt{T}_{n_1,r})
\ge \wt{T}_{n_1,r} - \wt{T}_{n_3,r}=\sum_{e\in\cR_u(n_1,n_3)}T_{e},
\ea
\edm
which concludes the proof.
\halmos
\endproof

\begin{proposition}
If $(\hat{x}^*,\hat{p}^*,\hat{y}^*,\hat{y}^{\prime*},\hat{w}^*)$ is an optimal solution of \eqref{opt:DT},
the output $\{\hat{\tau}^*_{u,e}:u\in\cV,\;e\in\cR_u\}$ of \algoref{alg:sol-DT-transform} with the input 
$(\hat{x}^*,\hat{p}^*,\hat{y}^*,\hat{y}^{\prime*},\hat{w}^*)$ is an optimal solution of the corresponding
tree-FRCVP-pause instance.
\end{proposition}
\proof{Proof.} 
Claim 1. The solution $\hat{\tau}^*$ to the tree-FRCVP-pause
instance gives the same objective value as the objective value of \eqref{opt:DT} yielded by the 
solution $(\hat{x}^*,\hat{p}^*,\hat{y}^*,\hat{y}^{\prime*},\hat{w}^*)$. Proof. It suffices to verify the following
properties: (a) For every $(e,t)$ pair, the departure time $\hat{\tau}^*_{u,e}$ is the same for
all $u\in\cV^*_{e,t}$, where $\cV^*_{e,t}:=\{u\in\cV:\;\hat{x}^*_{u,e,t}=1\}$. 
(b) For every $u\in\cV$, $\sum_{t\in\cS_u}\hat{p}^*_{u,t}-1$ equals the number of stops at intermediate nodes
along $\cR_u$ yielded by the solution $\hat{\tau}^*$. 

To verify (a), suppose $e=(n_1,n_2)$ and $t=[\ul{t},\ol{t}]$. 
The departure time from node $n_1$ is $\hat{\tau}^*_{u,e}=\ul{t}-\wt{T}_{n_1,r}$ for every $u\in\cV^*_{e,t}$,
which is the same for all vehicles in $\cV^*_{e,t}$. To verify (b), suppose $t_1\prec{t_2}\prec\ldots\prec{t_k}$ are 
time buckets from $\cS_u$ satisfying $\hat{p}_{u,t_i}=1$ for $i\in\{1,\ldots,k\}$ and $\hat{p}_{u,t}=0$ 
for any $t\in\cS_u\setminus\{1,\ldots,k\}$. 
Consider two consecutive time buckets $t_i:=[\ul{t}_i,\ol{t}_i]$ and $t_{i+1}:=[\ul{t}_{i+1},\ol{t}_{i+1}]$ in the above sequence.
There must exist two consecutive arcs $e_1\prec{e_2}\in\cR_u$ such that $e_1$  
satisfying $\hat{x}^*_{u,e_1,t_i}=\hat{x}^*_{u,e_2,t_{i+1}}=1$.
Suppose $e_1=(n_1,n_2)$ and $e_2=(n_2,n_3)$. The departure time of $u$ from node $n_1$ and $n_2$
are $\hat{\tau}^*_{u,n_1}=\ul{t}_i-\wt{T}_{n_1,r}$ and $\hat{\tau}^*_{u,n_2}=\ul{t}_{i+1}-\wt{T}_{n_2,r}$, respectively. 
Then the time difference is
\bdm
\hat{\tau}^*_{u,n_2}-\hat{\tau}^*_{u,n_1}=(\ul{t}_{i+1}-\wt{T}_{n_2,r}) - (\ul{t}_i-\wt{T}_{n_1,r})>\wt{T}_{n_1,r}-\wt{T}_{n_2,r}
=T_{e_1},
\edm
which shows that $u$ must stop at node $n_2$. Since any two consecutive time buckets from the sequence 
$t_1\prec{t_2}\prec\ldots\prec{t_k}$ induce exactly one stop, the whole sequence of time buckets induce
$\sum_{t\in\cS_u}\hat{p}^*_{u,t}-1$ stops according to $\hat{\tau}^*$. This concludes the proof of Claim 1.

Claim 2. Every feasible solution $\hat{\tau}$ of the tree-FRCVP-pause instance can be mapped to a feasible solution
$(\hat{x},\hat{p},\hat{y},\hat{y}^\prime,\hat{w})$ of \eqref{opt:DT} such that the objective value computed from $\hat{\tau}$ 
(using the tree-FRCVP-pause objective function)
is no greater than that computed from $(\hat{x},\hat{p},\hat{y},\hat{y}^\prime,\hat{w})$ (using the \eqref{opt:DT} objective function). 
Proof. For every $u\in\cR_u$ and $e=(n_1,n_2)\in\cR_u$, compute the relative departure time of $u$ from the starting node $n_1$ of $e$
as $\rtau_{u,e}=\hat{\tau}_{u,e}+\wt{T}_{n_1,r}$. Then exactly one of the following three cases holds:
\begin{itemize}
	\item[(a)] There is a unique time bucket $t=[\ul{t},\ol{t}]\in\cS$ such that $\rtau_{u,e}=\ul{t}=\ol{t}$.
	\item[(b)] $\hat{\tau}_{u,e}$ is contained in the first element $t=[\ul{t},\ol{t}]$ of $\cS_u$, i.e., $\rtau_{u,e}\in[\ul{t},\ol{t}]$.
	\item[(c)] There is a unique time bucket $t=[\ul{t},\ol{t}]\in\cS_u$ satisfying that $t$ is not the first element in $\cS_u$ 
and $\rtau_{u,e}\in(\ul{t},\ol{t}]$. 
\end{itemize}
In any case, we set $\hat{x}_{u,e,t}=1$. Compute values $\hat{p},\hat{y},\hat{y}^\prime,\hat{w}$ of the other variables
according to their definitions based on $\hat{x}$. Clearly, for any $e\in\cR$ and $u,v\in\cV_e$, if $(u,e)$ and $(v,e)$
are assigned to different time buckets, then $\rtau_{u,e}\neq\rtau_{v,e}$. This implies that the platoon formation
given by the solution $(\hat{x},\hat{p},\hat{y},\hat{y}^\prime,\hat{w})$ is no worse than that given by the solution $\hat{\tau}$.
For any $u$, consider any two consecutive arcs $e=(n_1,n_2)$ and $e^\prime=(n_2,n_3)$ from $\cR_u$.
If $\rtau_{u,e}=\rtau_{u,e^\prime}$, $(u,e)$ and $(u,e^\prime)$ must be assigned to the same time bucket,
and $u$ does not stop at $n_2$ either by the rule of $\hat{\tau}$ or by the rule of $(\hat{x},\hat{p})$.
Otherwise, if $\rtau_{u,e}<\rtau_{u,e^\prime}$, $u$ stops at $n_2$ by the rule of $\hat{\tau}$;
However, $(u,e)$ and $(u,e^\prime)$ could possibly be assigned to a same time bucket according to the
above construction of $\hat{x}$, and hence $u$ could possibly not stop at $n_2$ by the rule of $(\hat{x},\hat{p})$.
This indicates that the number of intermediate stops by the rule of $(\hat{x},\hat{p})$ is no more than that by
the rule of $\hat{\tau}$. Therefore, the Claim 2 holds. 

We now prove the proposition. Let $\hat{\tau}^*$ be an optimal solution of the tree-FRCVP-pause instance
with objective value being $obj^*$. 
Use Claim 2 to get a feasible solution $\hat{x}$ of \eqref{opt:DT} with objective value being $obj_1$. 
Let $\hat{x}^*$ be an optimal solution of \eqref{opt:DT}  with objective value being $obj_2$.
Use Claim 1 to get a feasible solution $\hat{\tau}$ of the tree-FRCVP-pause instance from $\hat{x}^*$
and let $obj_3$ be the objective value given by $\hat{\tau}$. The following relation holds:
\bdm
obj^*\overset{\textrm{Claim 2}}{\le}{obj_1}\le{obj_2}\overset{\textrm{Claim 1}}{=}{obj_3}\le{obj^*},
\edm
where the second inequality holds since $\hat{x}$ is a feasible solution 
but $\hat{x}^*$ is an optimal solution of \eqref{opt:DT}.
The last inequality holds since $\hat{\tau}$ is a feasible solution 
but $\hat{\tau}^*$ is an optimal solution of the tree-FRCVP-pause instance. 
This concludes the proof.
\halmos
\endproof

\begin{proposition}\label{prop:no-pause-interm}
There exists an optimal solution of the tree-FRCVP-pause problem such that
for every $v\in\cV$, a waiting decision is made at an intermediate node 
$n_2$ only if $\cV_{e_2}\neq\cV_{e_1}$, where 
$e_1=(n_1,n_2)$ and $e_2=(n_2,n_3)$ are two consecutive arcs on $\cR_v$. 
\end{proposition}
\proof{Proof.}
We only prove it for the case $\lambda=\infty$. 
The argument for the case $\lambda<\infty$ is similar. 
Consider an optimal solution $\tau=\{\tau_{u,e}: u\in\cV, e\in\cR_u\}$ of the instance,
where $\tau_{u,e}$ is the departure time of $u$ from the starting node of $e$.
Suppose there exists a vehicle $v$ such that a waiting option (with nonzero waiting time) 
is made at a node $n_2$ with $\cV_{e_2}=\cV_{e_1}$,
where $e_1=(n_1,n_2)$ and $e_2=(n_2,n_3)$ are two consecutive arcs on $\cR_v$.
Find the maximal path $P$ containing $e_1$ and $e_2$ on $\cR_v$ 
satisfying that $\cV_e=\cV_{e_1}$ for all $e\in{P}$. We introduce a notation
$\cV_P=\cV_{e_1}$ to indicate that the vehicle set is the same for all arcs on $P$.
We can divide $\cV_P$ into multiple classes based the following similarity 
relation: $u_1\sim{u_2}$ if $u_1$ and $u_2$ are in a same platoon in some $e\in{P}$.
Let $\cU$ be the class that contains $v$. Let $f_{\tau}(\cU,e)$ for $e\in{P}$ 
be the amount of fuel saving at arc $e$ achieved by platoon(s) within vehicles from $\cU$
under the solution $\tau$. Note that $f_{\tau}(\cU,e)$ is proportional to the length of $e$,
and hence we define $g_{\tau}(\cU,e)=f_{\tau}(\cU,e)/T_e$ to be the fuel-saving index 
which only depends on the platoon configuration of $\cU$ at $e$. Let
\bdm
e^*=\argmax_{e\in{P}}g_\tau(\cU,e),\quad g^*=\max_{e\in{P}}g_\tau(\cU,e).
\edm  
Suppose $e_0=(k_1,k_2)$ is the first arc on $P$ and $e^*=(k_3,k_4)$. Let
\bdm
\Delta_u = \tau_{u,e^*} - \tau_{u,e_0} - \sum_{e\in{P(k_1,k_3)}}T_e\quad\forall u\in\cU,
\edm
where $P(k_1,k_3)$ is the sub-path of $P$ connecting node $k_1$ and $k_3$.
Note that if $\Delta_u>0$, then $u$ has waited for at least once on $P$.  
We can construct a modified solution $\tau^\prime$
satisfying the following setting:
\bdm
\ba
&\tau^\prime_{u,e_0} = \tau_{u,e_0} - \Delta_u \quad\forall u\in\cU, \\
&\tau^\prime_{u,e^\prime} = \tau^\prime_{u,e}+T_e \qquad\forall u\in\cU, \forall \textrm{ consecutive }e,e^\prime\in{P},   
\ea
\edm
and all the other values of $\tau^\prime$ are set to be the same as $\tau$.
By construction, the solution $\tau^\prime$ is still feasible. 
The number of waiting of any vehicle in $\cU$ under $\tau^\prime$
is no greater than that under $\tau$, and all vehicles from $\cU$
do not wait at any intermediate node in $P$. By construction,
vehicles set $\cU$ can achieve the fuel-saving index $g^*$
at every arc of $P$ under $\tau^\prime$. This shows that
the solution $\tau^\prime$ yields an objective that is no-worse than $\tau$.
We can repeat the above procedure until for every maximal dir-path $P^\prime$ 
having the same vehicle set $\cV_{P^\prime}$ on it, 
vehicles in $\cV_{P^\prime}$ do not wait at any intermediate node of $P^\prime$.
Then we obtain a solution with desired property. 
\halmos
\endproof
\begin{remark}
\propref{prop:no-pause-interm} can help reducing the number of variables and constraints 
in \eqref{opt:DT} and \eqref{opt:VA}.
\end{remark}

\section{Numerical Investigation}
\label{sec:num-invest}
The proposed reformulations and algorithms are implemented using the Gurobi-Python API, 
and the effectiveness of \eqref{opt:VA} and \eqref{opt:DT} is tested against the MILP established 
in \cite{luo2020-route-then-sched-vehicle-platn} 
that is based on a continuous-time formulation with big-M coefficients involved 
(given in the appendix to be self-contained).

\subsection{Computational performance comparison}
\label{sec:comp-perf}
We generate 30 tree-FRCVP instances to test the computation performance of the proposed formulations
in this paper against the continuous-time based formulation (\ref{EC:cont-time-form}). 
Details of the instance generation process are provided in \ref{sec:prob-gen}.  
We compare the computational performance of using the 
assignment \eqref{opt:VA} and continuous-time \eqref{opt:cont_time} formulations 
to solve the numerical instances using a single core 2.3 GHz CPU with 32 GB memory.  
For each instance, the time limit of computing is one hour. 
We report the solution time (if the instance is solved within the time limit), the best objective value
and the relative optimality gap at 10 mins and 1 hour, respectively. The results  
are given in Table~\ref{tab:comp-perf}, where objective values are measured using the 
objective of the 100-vehicle system as a unit.
The general trend is that solving the instances with no constraint on the
platoon size is easier than those with capacity constraints for both formulations,
which is as expected. 

The results show that the assignment formulation outperforms the continuous-time
formulation in almost every numerical instance. 
In particular, using the assignment formulation, 
6 (out of 15) instances (Art100, Art200, CC100, CC200, CF100 and CF200) 
are solved to optimality for the non-capacitated case 
compared to the 4 instances (Art100, CC100, CC200 and CF100) for the capacitated case. 
While using the continuous-time formulation, only 2 instances (CC100, CF100) 
are solved to optimality for both the non-capacitated and capacitated cases.
Within the instances that can be solved to optimality using both formulations, 
a significant reduction in computational time is reported for the assignment formulation 
(at least 6 times faster than the continuous-time formulation), except for the instance (CC100, $\lambda=10$).
Within the instances that cannot be solved to optimality, the assignment formulation outperforms
the continuous-time formulation in the objective value and the optimality gap with a few
exceptions on instances with capacity constraints. At the 10 mins checkpoint, 
the assignment formulation wins in 17 instances, loses in 7 instances which are all in the capacitated case
(Art200, Art300, CF500, CC500, CF300, CF400 and CF500 with $\lambda=10$),
and the two formulations are tie in 6 instances. 
At the one-hour time limit, the assignment formulation loses in only 1 instance (CF500, $\lambda=10$).

We have conducted similar computational experiments to compare the performance
of the discrete-time formulation \eqref{opt:DT} and the continuous formulation \eqref{opt:cont_time}
modified to allow the stop-and-wait option on solving 30 numerical instances of the tree-FRCVP-pause
problem defined in \secref{sec:allow-waiting}. The results are given in \tabref{tab:pause-comp-perf}. 
It shows that \eqref{opt:DT} still has a certain advantage over \eqref{opt:cont_time} on solving
instances with the vehicle number $N\le250$. There are 19 problem instances solved to optimality
using formulation \eqref{opt:DT} compared to 12 problem instances using formulation \eqref{opt:cont_time}.
Within the 12 instances solved using both formulations, the formulation \eqref{opt:DT} achieves it 
with a shorter computational time, while for 11 instances that are not solved to optimality
within a 1-hour time limit, the formulation \eqref{opt:DT} either gives a better objective value
or a better optimality gap. Nevertheless, we have found that the formulation \eqref{opt:DT} 
becomes less effective when the vehicle number $N\ge300$, due to the number of variables
increases quickly with $N$. 

In the Gurobi solver, one needs to specify the method for solving the root relaxation
linear program, and the default is to use the dual simplex method. In most cases,
one can rely on the default method. But it turns out that in our case, the method selection 
plays a crucial role to make the two formulations perform at their best, especially for instances
involving a large number of vehicles (i.e., $N\ge300$), and the reason for which is interestingly related to 
the effectiveness of the formulation. The numerical investigation shows that the dual simplex method is 
sufficiently effective for solving the root LP using the continuous-time formulation,
while the barrier method (interior point method) is most effective for using the assignment formulation. 
(The results reported in Table~\ref{tab:comp-perf} are based on using the most effective LP-solving method for 
the two formulations.)
To get some insight, we have conducted an independent investigation of 
the three options of LP-solving methods: the primal simplex, dual simplex
and barrier methods on solving the root-relaxation linear programs of the two formulations.
The results are presented in Table~\ref{tab:LP-solve}. For the continuous-time
formulation, the primal and dual simplex methods are roughly on-par for solving 
the root LP's, whereas the barrier method takes longer time compared to the former
two methods, especially on CF300, CF400 and CF500. However, the outcome is very different
for the assignment formulation, for which case the barrier method remarkably 
outperforms the other two methods in all instances investigated, and there are 4$\sim$5
instances that can not be solved to optimality within the one-hour time limit by the primal and
dual simplex methods. The linear program corresponding to the assignment formulation seems to be much 
harder to solve compared to the continuous-time-formulation case. 
But once the assignment LP is solved, it yields a much more effective LP solution 
compared to the continuous-time LP. This is observed from that the objective ratio (defined in the caption of
Table~\ref{tab:LP-solve}) given by the assignment LP is 10\%$\sim$20\%
smaller than the one given by the continuous-time LP. Such a `no free lunch' observation 
is in fact a clear evidence that the assignment formulation is `tight and compact', and
is much stronger than the continuous-time formulation. Indeed, in the extreme situation,
the most effective linear program corresponds to the convex hull of all feasible solutions. 
But that polytope can have exponentially many extremal points (vertices),
in which case iterating from one vertex to another in its neighborhood (i.e., the approach used by the 
primal and dual simplex methods) can be extremely ineffecient for large instances. However, the high complexity 
of vertex enumeration on the surface can be
successfully avoided by the barrier method since it always iterates within the interior of the polytope.
Following this logic, if a formulation induces a polytope that is very close to the convex hull of feasible
points, it is very likely to have a similar computational behavior described above, which is exactly what
has been observed for the assignment formulation, showing that the assignment formulation is indeed
very strong.

\subsection{Fuel saving and time-window size correlation}
\label{sec:fuel-rtw}
\begin{figure}
\centering
\subfigure[]{\includegraphics[width=0.48\linewidth,trim=0cm 0.5cm 0cm 0cm]{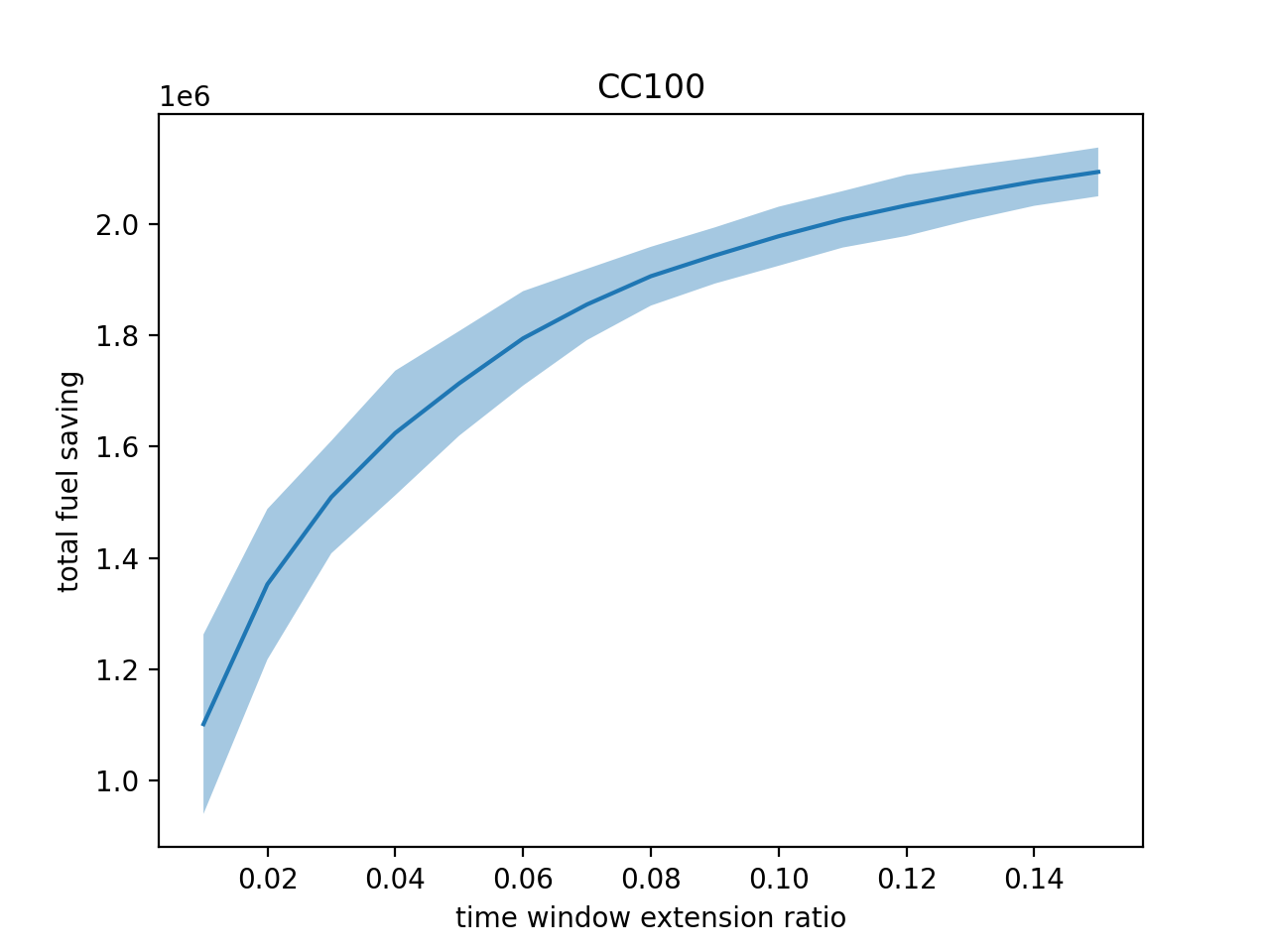}}
\hfill
\subfigure[]{\includegraphics[width=0.48\linewidth,trim=0cm 0.5cm 0cm 0cm]{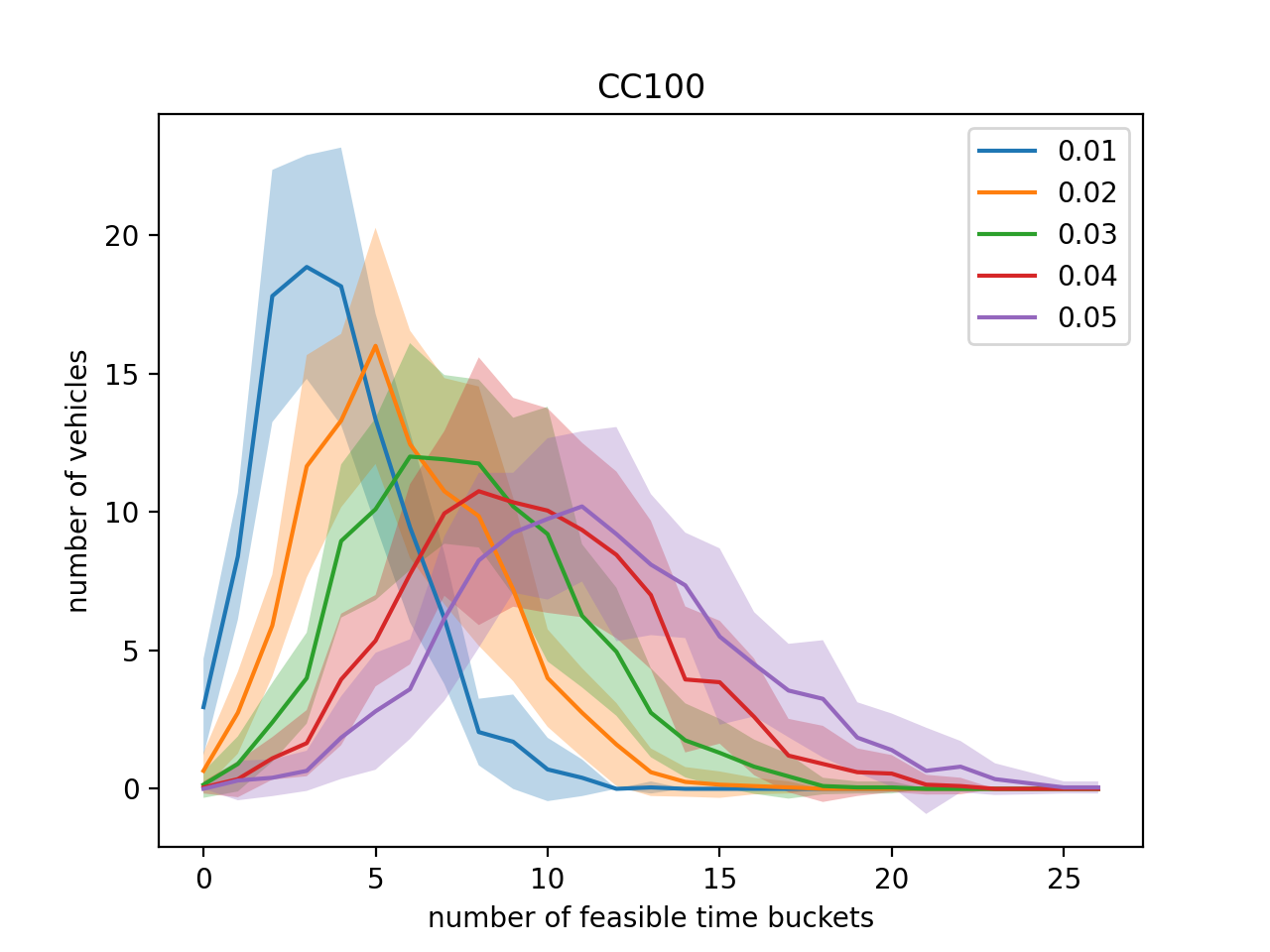}}
\captionsetup{justification=centering}
\caption{\scriptsize (a) The correlation between total fuel saving and the extension ratio ($\gamext/\gamfull$) on the instance CC100.
The extension ratio has been sequentially increased from 0.01 to 0.15 with the step size 0.01.
For a fixed extension ratio, 20 independent numerical instances are generated by randomly 
choosing $T^O_v$ from the interval $[0,\gamfull\ol{L}]$ for every vehicle. The sample mean
and standard deviation are computed for the objective obtained from the 20 samples.
The blue curve represents the mean value with respect to the extension ratio, and the shaded
region is corresponding to the $\pm 2\hat{\sigma}$ confidence interval, where $\sigma$ is the standard deviation.
(b) The distribution of vehicles over different values of the cardinality $|\cS_u|$ defined in Table~\ref{tab:var-TWOF}
(the number of time buckets that are feasible to a vehicle $u$ after running Algorithm~\ref{alg:ATD}). 
The plot is corresponding to the instance CC100. 
The distribution of vehicle number is plotted for the extension ratio
 $\gamext/\gamfull=0.01,\; 0.02,\; 0.03,\; 0.04 \text{ and } 0.05$ within the same figure. For each value of $\beta$, 
20 random instances are generated. The solid line represents the distribution of mean value
and the shaded region shows the $\pm 2\hat{\sigma}$ confidence interval, where $\hat{\sigma}$
is the sample standard deviation.}
\label{fig:fuel-vs-idle-time}
\end{figure}

From \ref{sec:EC-origin-dest-time}, the ratio $\gamext/\gamfull$ 
(defined as time-window-size / travel-time, referred as the extension ratio in this section) 
used in generating $T^O_v$ and $T^D_v$ for $v\in\cV$ measures the overall
time flexibility for platooning. 
We investigate the total fuel saving of the vehicle system as a function
of the extension ratio. In the experiments, we let the extension ratio grow from 0.01 to 0.15
with the step size 0.01. For a fixed extension ratio, we generate 20 numerical instances by 
randomly sampling $T^O_v$ for every vehicle $v$ and solve them to obtain the total fuel saving.
The experiments are conducted on the CC100 instance.
The results are given in \figref{fig:fuel-vs-idle-time}(a).
The figure shows that the total fuel saving monotonically increases as the extension ratio increases,
but the magnitude of increment decreases. It is also observed that the variation of the objective due to
the randomness of $T^O_v$ is less than 5\%. 
The variation of objective decreases as the extension ratio increases,
which shows that the increment of spare time can mitigate the uncertainty in the origin time selection.

We investigate the quantity $|\cS_v|$ (the number of feasible time buckets of a vehicle $v$)
versus the extension ratio, and the distribution of vehicles over the number of feasible time buckets. 
The results are given in \figref{fig:fuel-vs-idle-time}(b) for CC100. 
The distribution is generated by counting the number
of vehicles (within the 100-vehicle system) that have exactly $k$ feasible time buckets with $k$
from 1 to 30, and for the extension ratio being 0.01, 0.02, 0.03, 0.04 and 0.05, respectively. 
Intuitively, \figref{fig:fuel-vs-idle-time}(b) shows that the distribution increasingly 
spreads out as the extension ratio increases, meaning that more vehicles will
have a larger amount of feasible time buckets as the RTW size increases.    

\begin{figure}
\centering
\subfigure[]{\includegraphics[width=0.48\linewidth,trim=0cm 0.5cm 0cm 0cm]{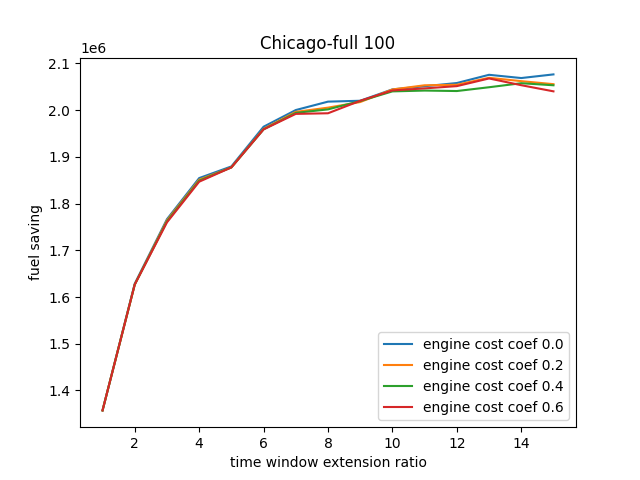}}
\hfill
\subfigure[]{\includegraphics[width=0.48\linewidth,trim=0cm 0.5cm 0cm 0cm]{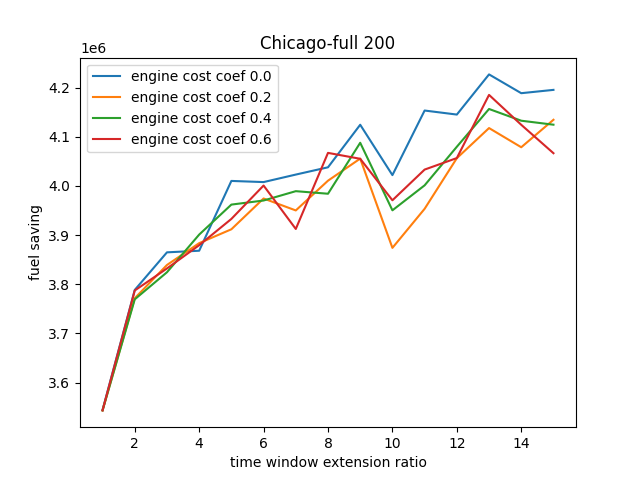}}
\captionsetup{justification=centering}
\caption{\scriptsize Fuel saving versus time window extension ratio under different coefficient of engine-restart 
fuel cost. The engine-restart fuel cost $\alpha_{\textrm{eng}}$ in \eqref{opt:DT} 
is set to be $c$ times the average fuel cost of passing an edge 
in the Chicago-full network, where $c\in\{0.0, 0.2, 0.4, 0.6\}$.
Each instance has been solved for 1hr with the formulation \eqref{opt:DT} using a one-core CPU.}
\label{fig:fuel-saving-diff-pause-cost}
\end{figure}

For the tree-FRCVP-pause problem, 
we have also studied the relation between the total fuel cost and time-window extension ratio
under different values of the engineer restart cost $\alpha_{\textrm{eng}}$. 
For a given $\alpha_{\textrm{eng}}$ (four options), we plot the total fuel saving obtained from solving
\eqref{opt:DT} for time-window extension ratio $\gamext/\gamfull$  
ranging from 1 to 15, with a 1-hour time limit on each problem instance.  
The results are shown in \figref{fig:fuel-saving-diff-pause-cost}. 
For the Chicago-full network with $N=100$ problem setting \figref{fig:fuel-saving-diff-pause-cost}(a),
each curve is almost monotonically increasing. But this is not the case for 
$N=200$ \figref{fig:fuel-saving-diff-pause-cost}(b). Theoretically, larger time-window size
should lead to higher fuel saving. But there could be a gap between what the solver can
find in one hour and the true optimal value and the gap depends on the complexity of the 
problem instance. In particular, \figref{fig:fuel-saving-diff-pause-cost}(b) shows that the 
high complexity problem instances are at $\gamext/\gamfull=7,8,10,14$, etc, and the complexity
is also $\alpha_{\textrm{eng}}$ dependent.

\section{Concluding Remarks}\label{sec:discussion}
The efficient time discretization algorithm and the novel vehicle-to-time-bucket
assignment formulation are proven to be more effective than the continuous-time 
formulation from the theoretical and computational perspectives, for tree-FRCVP instances. 
This approach could be extended to general instances to produce 
high-quality heuristic methods. 
The combination of this work and our previous one has concluded the framework 
of coordinated vehicle platooning over a road network with time constraints in the deterministic
setting. Further research could consider incorporating travel time uncertainty and the impact of traffic
under the framework of stochastic and robust optimization.

\section*{Acknowledgments}
This material is based upon work supported by the
U.S.\ Department of Energy, Office of Science, under contract number DE-AC02-06CH11357.
We are grateful for the data of greater Chicago highway network provided by Dr. Jeffrey Larson 
from the Argonne National Laboratory for generating numerical instances.

\section*{Appendix} 
\subsection*{A continuous-time formulation of the FRCVP problem}
\label{EC:cont-time-form}
The FRCVP problem has a continuous-time formulation developed in \cite{luo2020-route-then-sched-vehicle-platn}.
The formulation is given in \eqref{opt:cont_time} and the definition of notations used in the formulation is 
provided in Table~\ref{tab:cont_time}. Note that we assume the order of vehicles in a platoon
does not affect its fuel saving. In the continuous-time formulation, all vehicles are labeled with
positive integers and we only allow a vehicle with a greater index to follow a vehicle with a smaller
index as a convetion. (See the definition of the decision variable $f_{u,v,i,j}$ in Table~\ref{tab:cont_time}.)
\begin{subequations}
  \makeatletter
  \def\@currentlabel{CT}
  \makeatother
  \label{opt:cont_time}
  \renewcommand{\theequation}{CT.\arabic{equation}}
  \begin{alignat}{3}
    \text{maximize }     &\sum_{(i,j)\in\cup_v\cR_v} \Big(  \sum_{v\in\cV_{i,j}} \sigma^l C_{i,j} \ell_{v,i,j} + \sum_{u,v\in\cV_{i,j} u>v} \sigma^f C_{i,j} f_{u,v,i,j} \Big) \hspace{114pt}  \label{eqn:CT_1}
  \end{alignat}
  \vspace{-25pt}
  \begin{alignat}{3}
    \text{subject to: }  &  t_{v,O_v} \ge T_v^O    & \forall v\in\cV,  \label{eqn:CT_2} \\
& t_{v,D_v} \le T_v^D  & \forall  v\in\cV,   \label{eqn:CT_3}  \\
& t_{v,j} = t_{v,i} + T_{i,j} & \forall v\in\cV,\; \forall (i,j)\in\cR_v, \label{eqn:CT_4} \\
&  t_{u,i} - t_{v,i} \le M(1-f_{u,v,i,j}) & \forall (i,j)\in\bigcup_{v\in\cV} \cR_v, \; \forall u > v\in\cV_{i,j},  \label{eqn:CT_5}  \\
&  t_{u,i} - t_{v,i} \ge -M(1-f_{u,v,i,j}) & \hspace{1in} \forall (i,j)\in\bigcup_{v\in\cV} \cR_v, \; \forall u > v\in\cV_{i,j},  \label{eqn:CT_6}  \\
&\sum_{w\in\cV_{i,j}: w<v} f_{v,w,i,j} \le 1-\ell_{v,i,j} & \forall v\in \cV,\; (i,j)\in \cR_v,  \label{eqn:CT_7}  \\
&\sum_{u\in\cV_{i,j}: u>v} f_{u,v,i,j} \le (\lambda-1)\ell_{v,i,j} & \forall v\in \cV,\; (i,j)\in \cR_v,  \label{eqn:CT_8}  \\
&\sum_{u\in\cV_{i,j}: u>v} f_{u,v,i,j}\ge \ell_{v,i,j}  & \forall v\in \cV,\; (i,j)\in \cR_v,   \label{eqn:CT_9} \\
&f_{u,v,i,j}\in\{0,1\} & \forall u,v\in\cV, u>v, \forall (i,j)\in\cR_v,  \label{eqn:CT_11}  \\
&\ell_{v,i,j}\in\{0,1\} & \forall v\in\cV, \forall (i,j)\in\cR_v,  \label{eqn:CT_12} 
  \end{alignat}
\end{subequations}
where the $M$ in constraints \eqref{eqn:CT_5} and \eqref{eqn:CT_6} is a big-M coefficient 
satisfying $-M\le t_{u,i}-t_{v,i}\le M$ for all $(i,j)\in\cup_v\cR_v$, $u,v\in\cV_{i,j}$ with $u>v$ 
and any feasible values of $t_{u,i}$ and $t_{v,i}$.

The objective \eqref{eqn:CT_1} consists of two terms corresponding to 
fuel saving of lead and trailing vehicles on each arc, respectively.  
Note that maximizing the total fuel saving is equivalent to minimizing the total fuel cost. 
Constraints \eqref{eqn:CT_2}--\eqref{eqn:CT_4} ensure that vehicle travel times are consistent.
Constraints \eqref{eqn:CT_5}--\eqref{eqn:CT_6} ensure that if vehicle $u$ 
follows vehicle $v$ in a platoon on arc $(i,j)$, they must traverse the arc simultaneously.
Constraint~\eqref{eqn:CT_7} ensures that if vehicle $v$ is a lead vehicle 
on an arc, it cannot follow any other vehicles.
If $v$ is not a lead vehicle,
it can follow no more than one vehicle (the lead vehicle) on an arc. 
(Note that our formulation does not allow any trailing vehicle to lead another
vehicle.)
Constraint~\eqref{eqn:CT_8} enforces that
a vehicle can lead at most $\lambda-1$ other vehicles. 
Constraint~\eqref{eqn:CT_9} states that if $v$ is a lead vehicle
on an arc, the number of vehicles following $v$ on that arc must be at least one.
\begin{table}
  \footnotesize
  \caption{Sets, parameters, and variables for the continuous-time formulation. \label{tab:cont_time}}
  \begin{tabularx}{\textwidth}{lX}
    \toprule
    Set & Definition\\
    \midrule
			$\cR_v$ & route of vehicle $v\in\cV$ (an ordered subset of arcs in the route)\\
			$\cN_v$ & set of nodes on $\cR_v$, $v\in\cV$\\
      $\cV_{i,j}$ & set of vehicles taking the arc $(i,j) \in \cup_{v\in\cV}\cR_v$ \\
    \midrule\\[-.35cm]
    Parameter & Definition\\
    \midrule
    $T_v^O,T_v^D,T_{i,j}$ & defined in Table~\ref{tab:param1}\\
			$\lambda$ & maximum number of vehicles permitted in a platoon  \\
    $M_{u,v,i,j}$			& a vehicle-pair and arc dependent big-M coefficient \\
    \midrule\\[-.35cm]
    Variable & Definition\\
			$t_{v,O_v}$ & departure time of vehicle $v$ from node $O_v$ \\
			$t_{v,i}$ & arrival time of vehicle $v$ at an intermediate node $i$ along its route $\cR_v$ \\
      $f_{u,v,i,j}$ & 1 if vehicle $u$ follows vehicle $v$ at arc $(i,j)\in\cR_u\cap\cR_v$, 0 otherwise. (We require $u>v$ by the convention of ordering.) \\
			$\ell_{v,i,j}$ & 1 if vehicle $v$ leads a platoon at arc $(i,j)\in\cR_v$, 0 otherwise\\
    \bottomrule
  \end{tabularx}
\end{table}

\subsection*{Network structure}
\begin{figure}
\centering
\subfigure[ArtifNet (Art) network]{\includegraphics[width=0.45\linewidth,trim=1cm 2cm 1cm 0.5cm]{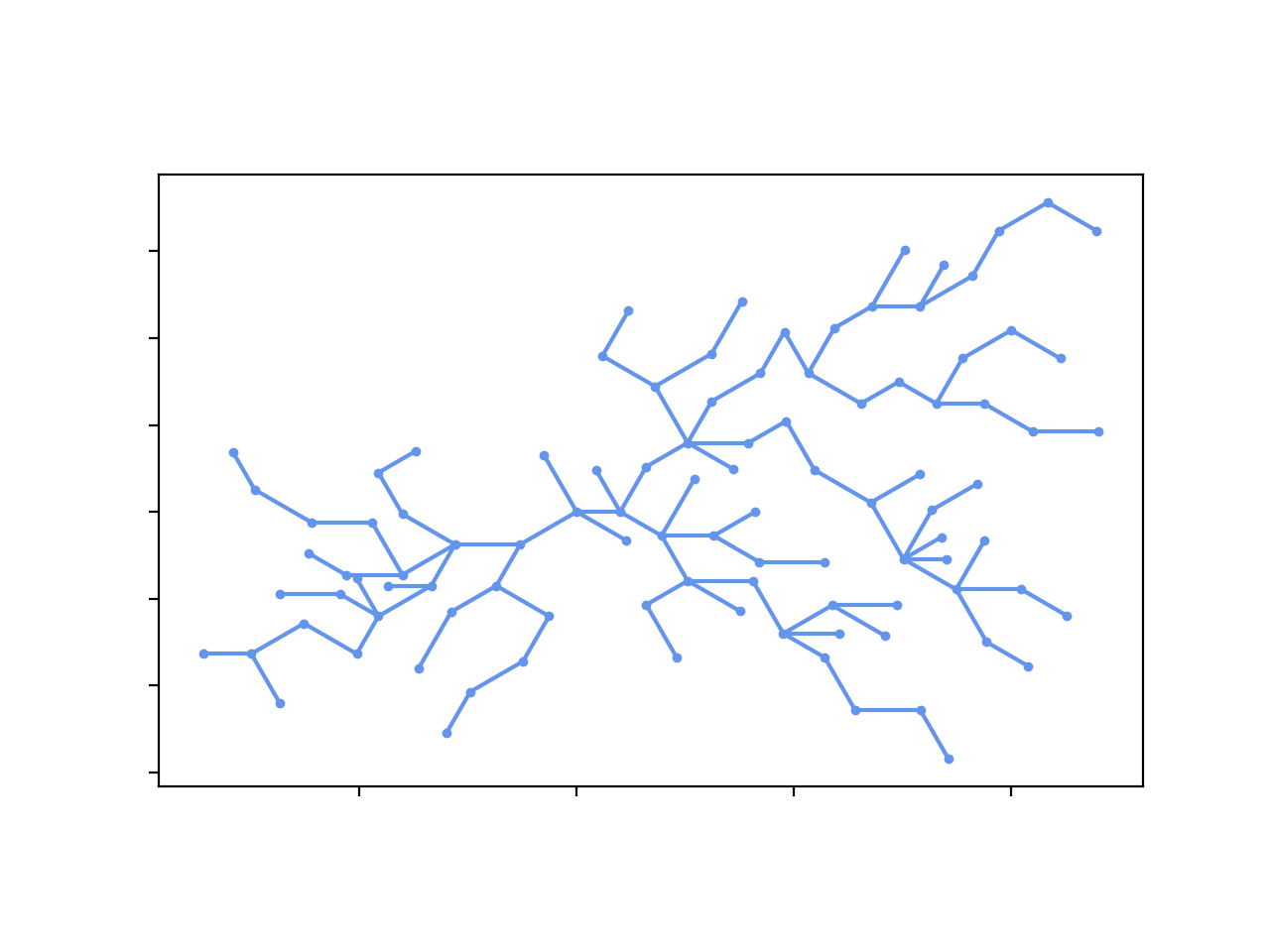}}
\hfill 
\subfigure[ArtifNet with 100 vehicles]{\includegraphics[width=0.45\linewidth,trim=1cm 2cm 1cm 0.5cm]{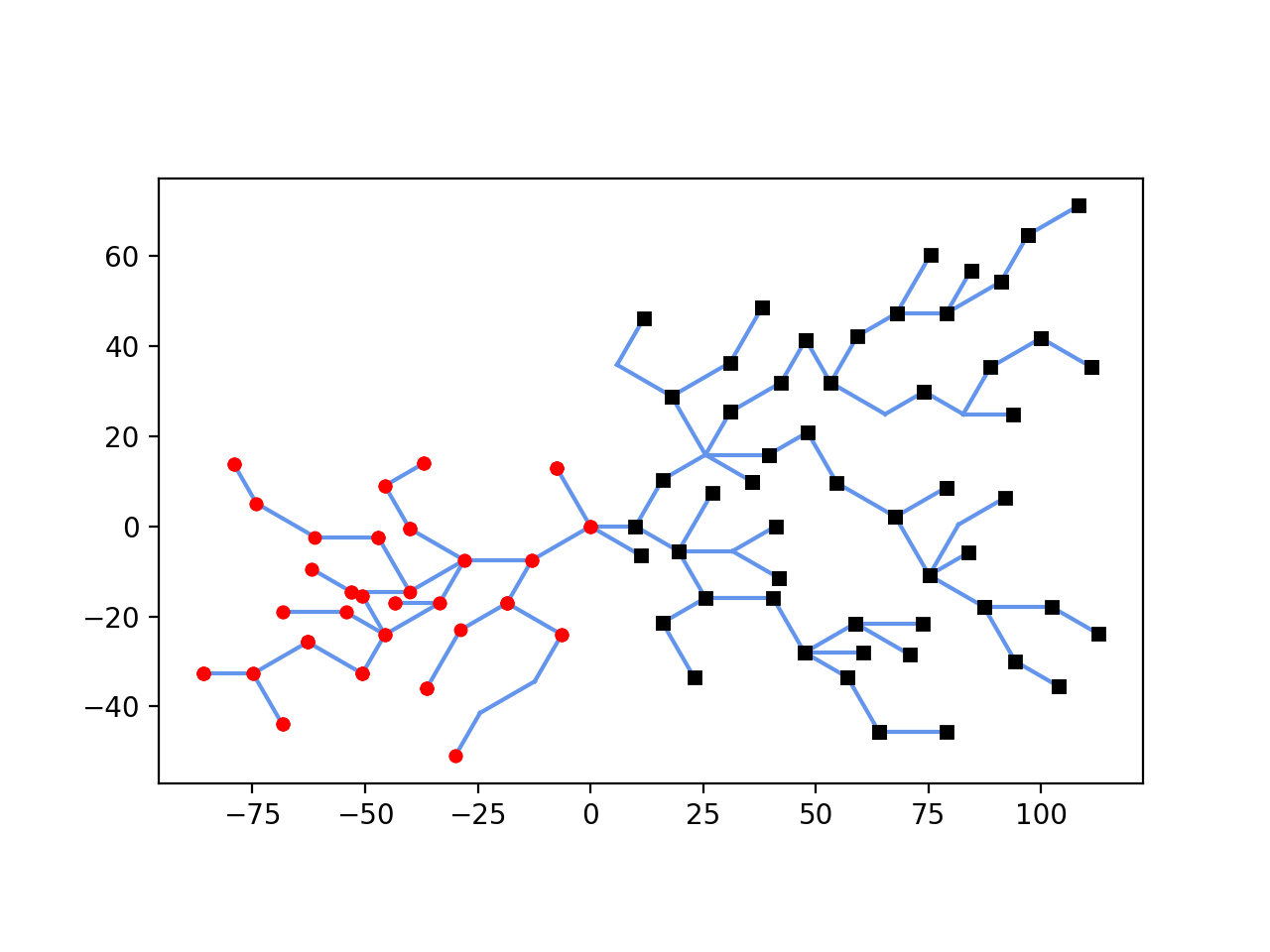}} 
\hfill
\subfigure[Chicago-cluster (CC) network]{\includegraphics[width=0.45\linewidth,trim=1cm 1.5cm 1cm 0.5cm]{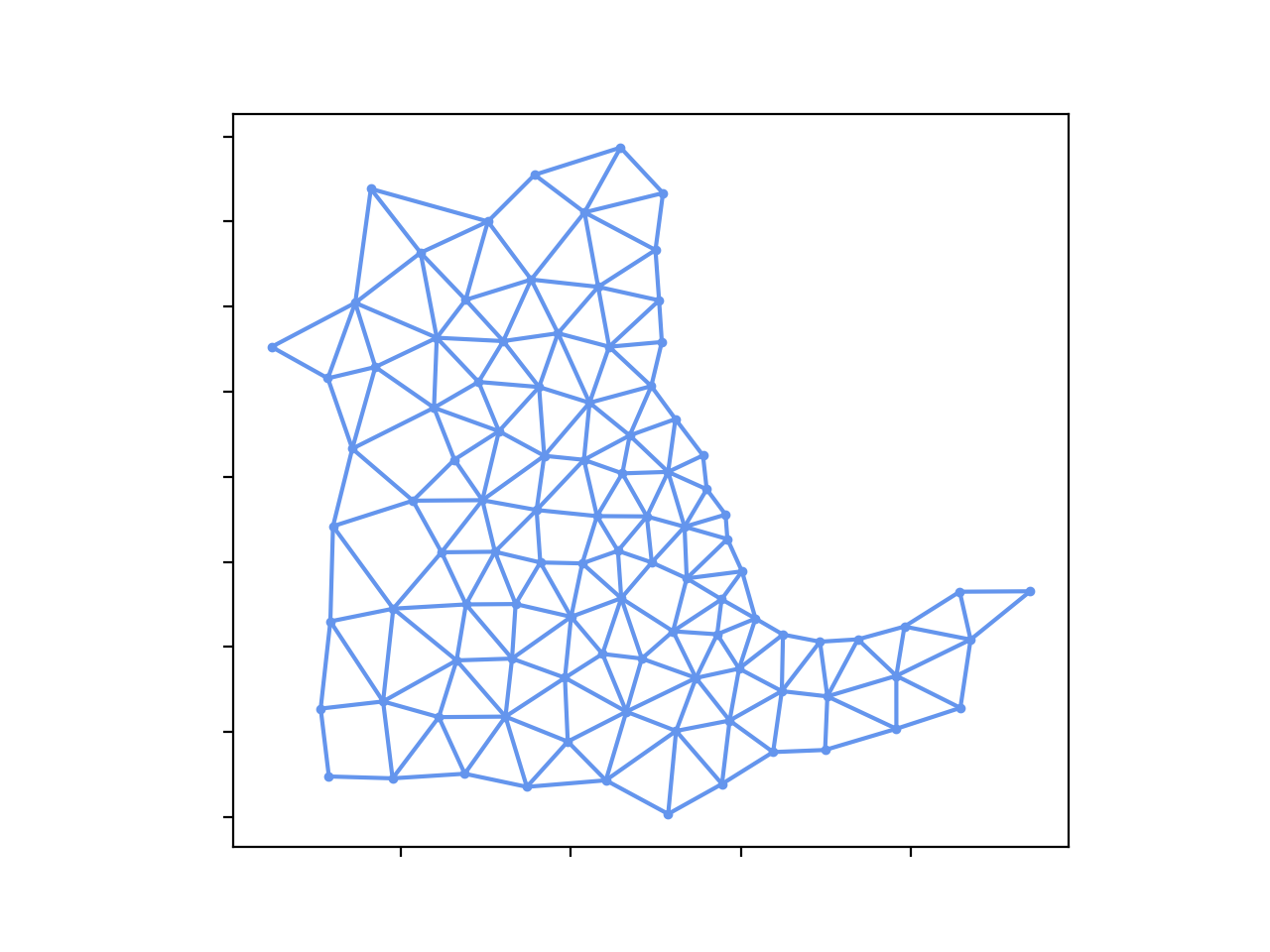}}
\hfill 
\subfigure[Chicago-cluster with 100 vehicles]{\includegraphics[width=0.45\linewidth,trim=1cm 1.5cm 1cm 0.5cm]{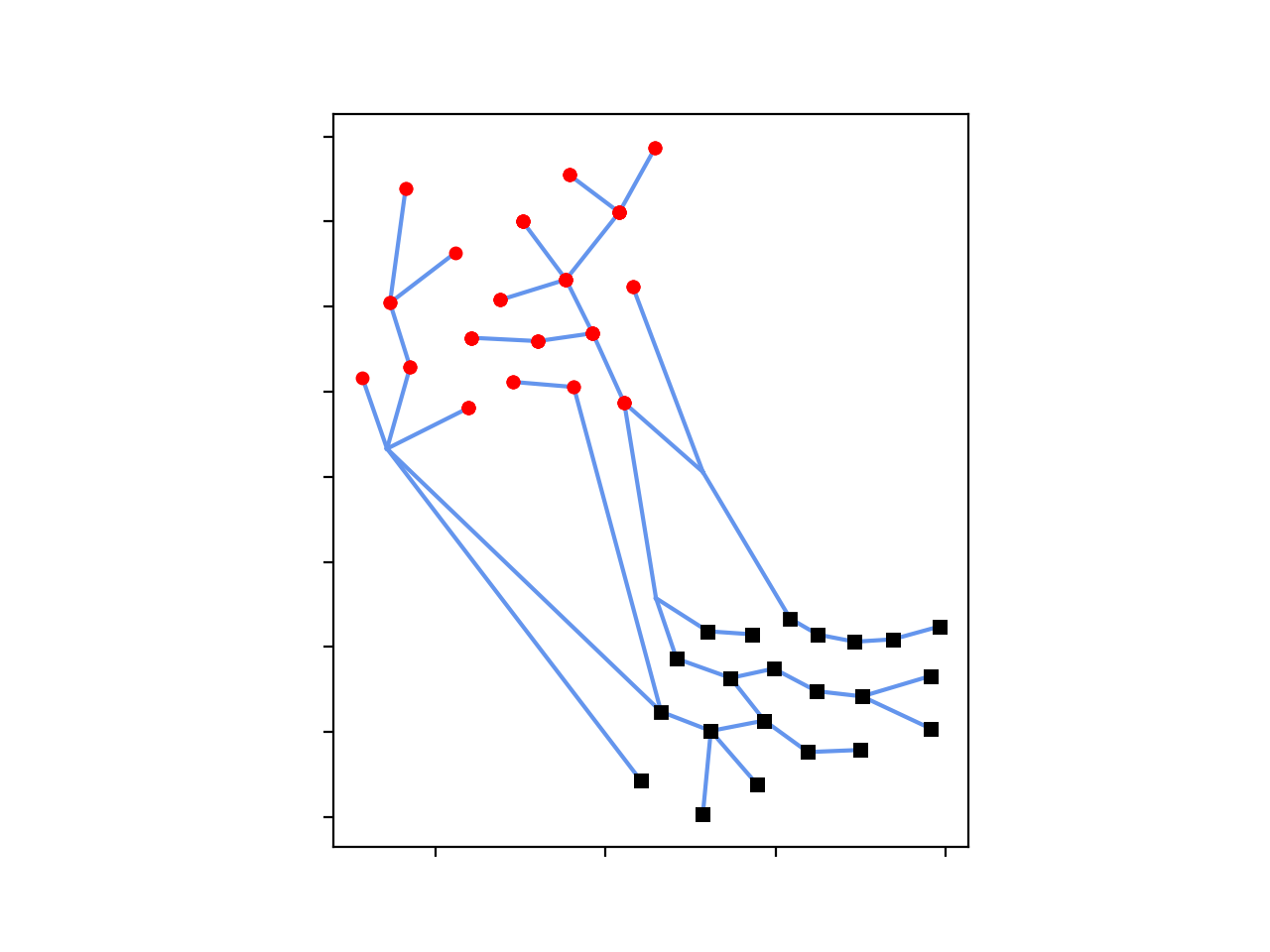}}
\hfill
\subfigure[Chicago-full (CF) network]{\includegraphics[width=0.45\linewidth,trim=1cm 1.5cm 1cm 0.5cm]{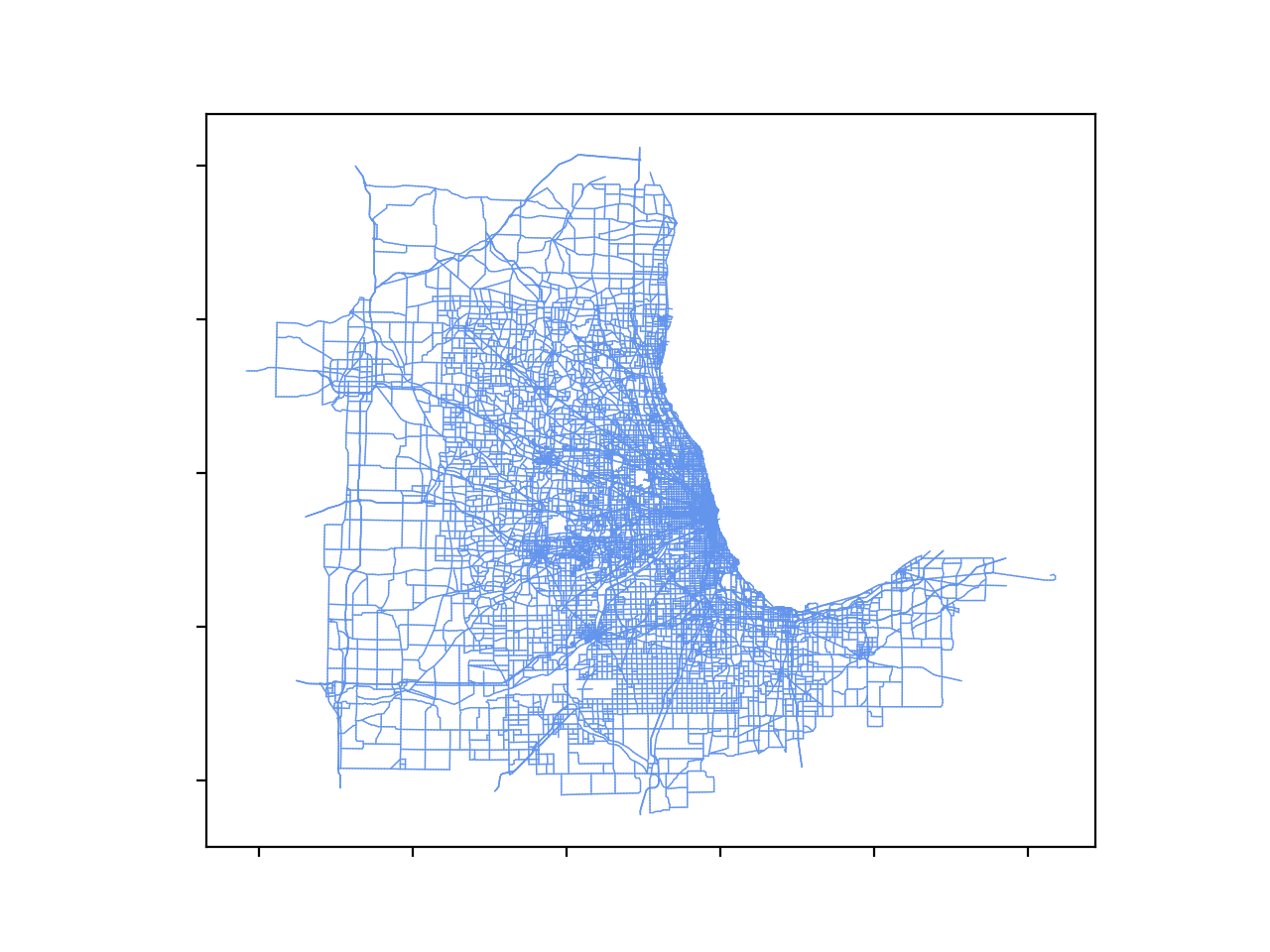}}
\hfill
\subfigure[Chicago-full with 100 vehicles]{\includegraphics[width=0.45\linewidth,trim=1cm 1.5cm 1cm 0.5cm]{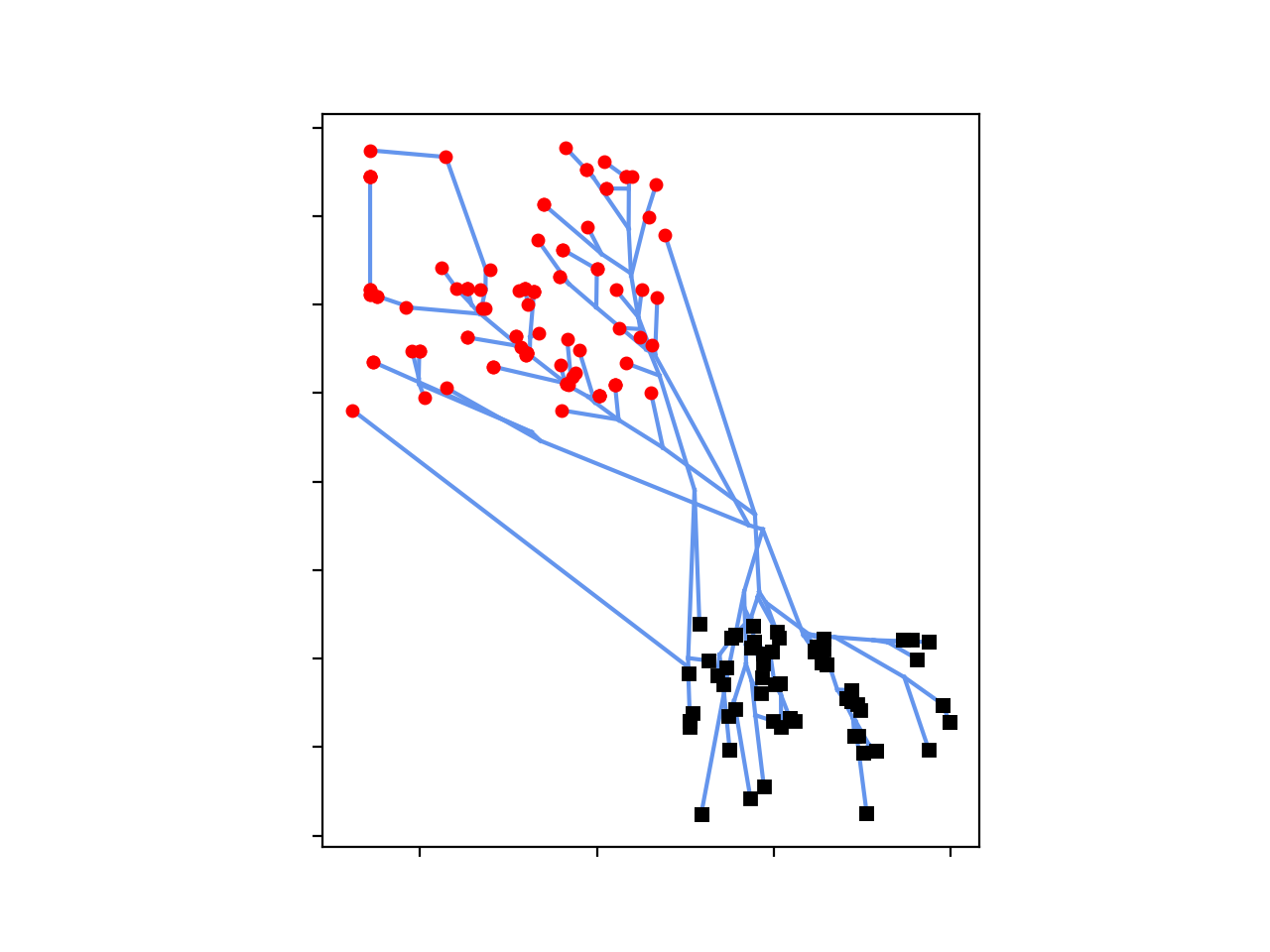}}
\captionsetup{justification=centering}
\caption{\scriptsize Part (a), (c) and (e) are networks used for generating numerical instances.
Part (b), (d) and (e) are examples of a 100-vehicle system defined on the three networks,
respectively, with circle and square dots representing the origin and destination nodes, respectively.
The original route of each vehicle is the shortest path, but minor adjustment on routes has been made to 
ensure that the route-graph is a tree. Straight lines are used to link the origin and destination nodes
in (d) and (f) to simplify the plot.}\label{fig:network}
\end{figure}

The road networks generated in this work is given in \figref{fig:network}.

\begin{table}
\caption{\scriptsize Comparison of computational performance between the assignment formulation
and the continuous-time formulation \eqref{opt:cont_time} based on 30 numerical instances of FRCVP with the route graph being a polytree. 
The time limit for solving each instance is one hour
and the tolerance threshold for optimality gap is set to be $0.1\%$.
The solution time is reported if the instance can be solved within the one-hour time limit,
otherwise it is labeled as `-'. For the two formulations, the objective value and relative optimality gap (in percent) 
are reported after computing for 10 mins, and one hour, respectively. Note that the objective value
has been rescaled to increase the interpretability. Specifically, for a each road network,  
the optimal objective value of the 100-vehicle system is set as the reference value (a unit)
to represent the objective values of other instances defined on the same network. 
For example, the optimal objective of Art100-inf is set to be 1.0 as the reference for 
representing objective values of Art200-inf to Art500-inf. 
} \label{tab:comp-perf}
\begin{center}
{\scriptsize
\begin{tabular}{cc|ccccc|ccccc}
\hline\hline
   &  & \multicolumn{5}{c|}{Assignment Formulation} & \multicolumn{5}{c}{Continuous-time Formulation} \\
\cline{1-12}
   &  &  & \multicolumn{2}{|c}{10 min} & \multicolumn{2}{|c|}{1 hour} &  & \multicolumn{2}{|c|}{10 min} & \multicolumn{2}{c}{1 hour} \\
   \cline{3-12}
Instance &  Capacity & sol. time (s) & obj. & gap (\%) & obj. & gap (\%) & sol. time (s) & obj. & gap (\%) & obj. & gap (\%) \\
\hline
Art100	&	inf	&	14	&	1.000	&	0.00	&	1.000	&	0.00	&	 -	&	1.000	&	0.80	&	1.000	&	0.23	\\
Art200	&	inf	&	2454	&	2.169	&	0.25	&	2.170	&	0.00	&	 -	&	2.082	&	6.90	&	2.143	&	3.21	\\
Art300	&	inf	&	 -	&	3.496	&	1.77	&	3.523	&	0.70	&	 -	&	3.371	&	8.16	&	3.439	&	5.09	\\
Art400	&	inf	&	 -	&	4.631	&	4.29	&	4.778	&	0.97	&	 -	&	4.539	&	9.19	&	4.568	&	7.26	\\
Art500	&	inf	&	 -	&	5.891	&	5.10	&	5.972	&	3.66	&	-	&	5.802	&	9.64	&	5.838	&	7.97	\\
CC100	&	inf	&	4	&	1.000	&	0.00	&	1.000	&	0.00	&	35	&	1.000	&	0.00	&	1.000	&	0.00	\\
CC200	&	inf	&	155	&	2.371	&	0.00	&	2.371	&	0.00	&	 -	&	2.365	&	1.02	&	2.370	&	0.70	\\
CC300	&	inf	&	 -	&	3.818	&	0.23	&	3.819	&	0.05	&	 -	&	3.790	&	2.13	&	3.806	&	1.10	\\
CC400	&	inf	&	 -	&	5.253	&	0.29	&	5.253	&	0.18	&	 -	&	5.076	&	5.33	&	5.186	&	2.05	\\
CC500	&	inf	&	 -	&	6.683	&	1.02	&	6.719	&	0.31	&	 -	&	6.516	&	5.29	&	6.659	&	2.17	\\
CF100	&	inf	&	10	&	1.000	&	0.00	&	1.000	&	0.00	&	450	&	1.000	&	0.00	&	1.000	&	0.00	\\
CF200	&	inf	&	1291	&	2.328	&	0.12	&	2.329	&	0.01	&	 -	&	2.295	&	2.62	&	2.319	&	1.17	\\
CF300	&	inf	&	 -	&	3.680	&	1.71	&	3.721	&	0.44	&	 -	&	3.590	&	6.39	&	3.616	&	4.86	\\
CF400	&	inf	&	 -	&	5.129	&	0.50	&	5.132	&	0.20	&	 -	&	4.917	&	7.49	&	4.944	&	5.99	\\
CF500	&	inf	&	 -	&	6.320	&	2.64	&	6.400	&	0.47	&	 -	&	6.203	&	7.03	&	6.219	&	6.42	\\
\hline
Art100	&	10	&	407	&	1.000	&	0.10	&	1.000	&	0.10	&	 -	&	1.000	&	1.00	&	1.000	&	0.23	\\
Art200	&	10	&	 -	&	1.979	&	11.75	&	2.143	&	1.20	&	 -	&	2.032	&	9.89	&	2.140	&	4.01	\\
Art300	&	10	&	 -	&	3.144	&	13.60	&	3.438	&	3.32	&	 -	&	3.303	&	10.92	&	3.323	&	8.32	\\
Art400	&	10	&	 -	&	4.063	&	18.02	&	4.534	&	7.16	&	 -	&	na	&	na	&	4.470	&	9.62	\\
Art500	&	10	&	 -	&	4.948	&	22.00	&	5.733	&	8.41	&	 -	&	na	&	na	&	5.730	&	9.37	\\
CC100	&	10	&	118	&	1.000	&	0.00	&	1.000	&	0.00	&	74	&	1.000	&	0.00	&	1.000	&	0.00	\\
CC200	&	10	&	1228	&	2.357	&	1.62	&	2.370	&	0.00	&	 -	&	2.357	&	1.44	&	2.369	&	0.74	\\
CC300	&	10	&	 -	&	3.722	&	4.81	&	3.792	&	1.15	&	 -	&	3.701	&	4.44	&	3.782	&	1.79	\\
CC400	&	10	&	 -	&	5.000	&	6.76	&	5.131	&	3.22	&	 -	&	4.890	&	9.09	&	5.100	&	3.97	\\
CC500	&	10	&	 -	&	6.293	&	8.51	&	6.568	&	3.41	&	 -	&	6.385	&	7.07	&	6.472	&	4.94	\\
CF100	&	10	&	142	&	1.000	&	0.00	&	1.000	&	0.00	&	859	&	1.000	&	0.35	&	1.000	&	0.01	\\
CF200	&	10	&	 -	&	2.246	&	5.03	&	2.315	&	0.41	&	 -	&	2.241	&	5.66	&	2.313	&	1.60	\\
CF300	&	10	&	 -	&	3.276	&	14.35	&	3.619	&	3.58	&	 -	&	3.502	&	8.51	&	3.523	&	6.96	\\
CF400	&	10	&	 -	&	4.394	&	16.31	&	4.856	&	6.11	&	 -	&	4.783	&	8.86	&	4.821	&	7.02	\\
CF500	&	10	&	 -	&	5.695	&	13.58	&	5.948	&	8.78	&	 -	&	5.985	&	9.28	&	6.028	&	8.15	\\
\hline
\end{tabular}
}
\end{center}
\end{table}

\begin{table}
\caption{\scriptsize The solution time comparison of solving the root relaxation linear programs (LP's) 
using three different methods: the primal simplex, the dual simplex and the barrier method (interior point method)
for the assignment and continuous-time formulations, respectively, tested on numerical instances
involving large number of vehicles with $\lambda=\infty$. 
The time limit for solving the LP relaxation is one hour. 
If an instance cannot be solved to optimality by a method, the corresponding cell of the computational time is marked as `-'.
The column `obj. ratio' is the ratio between the optimal objective of the root relaxation LP and the optimal objective of
the instance. As a maximization problem, the ratio is always greater than or equal to 1.0.}\label{tab:LP-solve}
\begin{center}
{\scriptsize
 \begin{tabular}{c|rrrr|rrrr}
  \hline\hline
   		 & \multicolumn{4}{c|}{Assignment Formulation \eqref{opt:VA}} & \multicolumn{4}{c}{Continuous Time Formulation \eqref{opt:cont_time}} \\
 \cline{2-9}
 Instance	&	primal (sec.)	&	dual (sec.)	&	barrier (sec.)	&	obj. ratio	&	primal (sec.)	&	dual (sec.)	&	barrier (sec.)	&	obj. ratio	\\
 \cline{1-9}
Art300	&	1617	&	707	&	40	&	1.040	&	10	&	8	&	22	&	1.166	\\
Art400	&	-	&	-	&	91	&	1.035	&	54	&	25	&	47	&	1.143	\\
Art500	&	-	&	-	&	200	&	1.055	&	47	&	64	&	86	&	1.150	\\
CC300	&	239	&	66	&	10	&	1.024	&	3	&	4	&	11	&	1.144	\\
CC400	&	963	&	302	&	22	&	1.016	&	7	&	15	&	28	&	1.111	\\
CC500	&	-	&	1147	&	43	&	1.012	&	11	&	11	&	40	&	1.089	\\
CF300	&	2913	&	3392	&	96	&	1.060	&	20	&	30	&	67	&	1.200	\\
CF400	&	-	&	-	&	350	&	1.044	&	29	&	41	&	115	&	1.157	\\
CF500	&	-	&	-	&	499	&	1.049	&	75	&	83	&	215	&	1.155	\\
\hline
 \end{tabular}
}
\end{center}
\end{table}

\begin{table}
\caption{\scriptsize Comparison of computational performance between the discrete-time formulation
\eqref{opt:DT} and the continuous-time formulation \eqref{opt:cont_time} with the stop-and-wait option 
at intermediate nodes based on 30 numerical instances of the tree-FRCVP-pause problem. 
} \label{tab:pause-comp-perf}
\begin{center}
{\scriptsize
\begin{tabular}{cc|ccccc|ccccc}
\hline\hline
   &  & \multicolumn{5}{c|}{Discrete-time Formulation \eqref{opt:DT}} & \multicolumn{5}{c}{Continuous-time Formulation \eqref{opt:cont_time}} \\
\cline{1-12}
   &  &  & \multicolumn{2}{|c}{10 min} & \multicolumn{2}{|c|}{1 hour} &  & \multicolumn{2}{|c|}{10 min} & \multicolumn{2}{c}{1 hour} \\
   \cline{3-12}
Instance &  Capacity & sol. time (s) & obj. & gap (\%) & obj. & gap (\%) & sol. time (s) & obj. & gap (\%) & obj. & gap (\%) \\
\hline
Art50	&	inf	&	5	&	0.335	&	0.00	&	0.335	&	0.00	&	3	&	0.335	&	0.00	&	0.335	&	0.00	\\
Art100	&	inf	&	355	&	1.000	&	0.01	&	1.000	&	0.01	&	 -	&	1.000	&	0.84	&	1.000	&	0.51	\\
Art150	&	inf	&	1083	&	0.997	&	0.17	&	0.997	&	0.01	&	 -	&	0.985	&	2.64	&	0.985	&	2.43	\\
Art200	&	inf	&	 -	&	1.977	&	19.89	&	2.151	&	0.24	&	 -	&	2.028	&	9.05	&	2.123	&	3.38	\\
Art250	&	inf	&	 -	&	1.683	&	22.42	&	1.768	&	6.24	&	 -	&	1.733	&	6.45	&	1.754	&	4.65	\\
CC50	&	inf	&	2	&	0.481	&	0.00	&	0.481	&	0.00	&	1	&	0.481	&	0.00	&	0.481	&	0.00	\\
CC100	&	inf	&	17	&	1.000	&	0.00	&	1.000	&	0.00	&	147	&	1.000	&	0.01	&	1.000	&	0.01	\\
CC150	&	inf	&	101	&	1.652	&	0.00	&	1.652	&	0.00	&	669	&	1.652	&	0.11	&	1.652	&	0.00	\\
CC200	&	inf	&	 -	&	2.274	&	14.61	&	2.372	&	0.26	&	 -	&	2.368	&	1.08	&	2.370	&	0.84	\\
CC250	&	inf	&	 -	&	2.795	&	22.83	&	3.014	&	8.26	&	 -	&	3.056	&	2.15	&	3.074	&	1.09	\\
CF50	&	inf	&	3	&	0.430	&	0.00	&	0.430	&	0.00	&	3	&	0.430	&	0.00	&	0.430	&	0.00	\\
CF100	&	inf	&	48	&	1.000	&	0.00	&	1.000	&	0.00	&	829	&	1.000	&	0.21	&	1.000	&	0.01	\\
CF150	&	inf	&	821	&	1.611	&	1.03	&	1.625	&	0.00	&	 -	&	1.619	&	1.06	&	1.625	&	0.39	\\
CF200	&	inf	&	1730	&	2.180	&	19.43	&	2.328	&	0.01	&	 -	&	2.269	&	4.11	&	2.319	&	1.18	\\
CF250	&	inf	&	2524	&	2.828	&	16.82	&	2.954	&	0.00	&	 -	&	2.795	&	6.70	&	2.942	&	0.96	\\
\hline																							
Art50	&	10	&	6	&	0.335	&	0.00	&	0.335	&	0.00	&	3	&	0.335	&	0.00	&	0.335	&	0.00	\\
Art100	&	10	&	138	&	1.000	&	0.00	&	1.000	&	0.00	&	 -	&	0.999	&	0.79	&	1.000	&	0.43	\\
Art150	&	10	&	 -	&	0.997	&	0.33	&	0.997	&	0.11	&	 -	&	0.986	&	2.59	&	0.993	&	1.77	\\
Art200	&	10	&	 -	&	2.102	&	4.33	&	2.139	&	0.84	&	 -	&	1.989	&	11.33	&	2.105	&	4.87	\\
Art250	&	10	&	 -	&	1.758	&	5.74	&	1.780	&	3.02	&	 -	&	1.690	&	8.85	&	1.755	&	4.73	\\
CC50	&	10	&	3	&	0.481	&	0.00	&	0.481	&	0.00	&	1	&	0.481	&	0.00	&	0.481	&	0.00	\\
CC100	&	10	&	33	&	1.000	&	0.00	&	1.000	&	0.00	&	78	&	1.000	&	0.01	&	1.000	&	0.01	\\
CC150	&	10	&	689	&	1.652	&	0.01	&	1.652	&	0.01	&	1114	&	1.651	&	0.24	&	1.652	&	0.00	\\
CC200	&	10	&	 -	&	2.314	&	5.04	&	2.372	&	0.35	&	 -	&	2.367	&	1.43	&	2.371	&	0.84	\\
CC250	&	10	&	 -	&	2.966	&	6.41	&	3.070	&	1.25	&	 -	&	3.040	&	2.69	&	3.078	&	0.92	\\
CF50	&	10	&	4	&	0.430	&	0.00	&	0.430	&	0.00	&	2	&	0.430	&	0.00	&	0.430	&	0.00	\\
CF100	&	10	&	217	&	1.000	&	0.00	&	1.000	&	0.00	&	957	&	1.000	&	0.25	&	1.000	&	0.01	\\
CF150	&	10	&	2185	&	1.619	&	1.20	&	1.625	&	0.01	&	 -	&	1.621	&	0.91	&	1.622	&	0.58	\\
CF200	&	10	&	 -	&	2.292	&	3.28	&	2.324	&	0.19	&	 -	&	2.277	&	3.70	&	2.312	&	1.59	\\
CF250	&	10	&	 -	&	2.925	&	2.51	&	2.942	&	0.36	&	 -	&	2.787	&	6.88	&	2.926	&	1.83	\\
\hline
\end{tabular}
}
\end{center}
\end{table}

\def\bibfont{\footnotesize}
\bibliographystyle{informs2014} 
\bibliography{references}


\ECSwitch
\ECHead{Supplemental Material}

\section{Omitted Proofs}\label{sec:suppl-1}
\thmI*
\proof{Proof.}
Denote $T$ as the random variable of the first end point of a RTW which follows
the distribution $F_T$, and $G$ as the RTW length random variable
following the distribution $F_G$. 
Consider the RTW of a vehicle $u$, denoted as $[T_u, T_u + G_u]$. Suppose
we first put $\rtw_u$ into the coordinate axis, and  
put the RTW's of other vehicles into the axis one by one.
When a new $\rtw_v$ is added to the axis, it may intersect 
with $\rtw_u$ and if an end point of $\rtw_v$ drops into the range
of $\rtw_u$, it almost surely splits a current feasible time bucket into two, and hence
refining the set of feasible time buckets of $u$. Quantitively, if one (two) end point(s)
of $\rtw_v$ is (are) within the interval $\rtw_u$, the number of feasible time buckets of $u$
is increased by one (two). So it amounts to compute the probability of these events. 
  
Consider putting $\rtw_v$ (of vehicle $v$) into the axis and let $\Delta$ denote
the increment of feasible time buckets of $u$ right after the refinement of adding $\rtw_v$.
For the events $\{T_u<T_v<T_u+G_u<T_v+G_v\}$ and $\{T_v<T_u<T_v+G_v<T_u+G_u\}$,
the $\Delta=1$, while for the event $\{T_u<T_v<T_v+G_v<T_u+G_u\}$, $\Delta=2$.
The probability of the first event can be computed as
\bdm
\ba
&P(T_u<T_v<T_u+G_u<T_v+G_v)=P(T_1<T_2, T_2-T_1<G_1, G_1-G_2<T_2-T_1) \\
&=P(T_2-T_1<G_1, G_1-G_2<T_2-T_1|T_1<T_2)\cdot P(T_1<T_2) \\
&=\frac{1}{2}P(G_1>|T_1-T_2|, G_1-G_2<|T_1-T_2||T_1<T_2) \\
&=\frac{1}{2}P(G_1>|T_1-T_2|, G_1-G_2<|T_1-T_2|, G_1<G_2 |T_1<T_2) \\
&\quad + \frac{1}{2}P(G_1>|T_1-T_2|, G_1-G_2<|T_1-T_2|, G_1>G_2 |T_1<T_2) \\
&=\frac{1}{2}P(G_2>G_1>|T_1-T_2| |T_1<T_2) \\
&\quad + \frac{1}{2}P(G_1>|T_1-T_2|, |G_1-G_2|<|T_1-T_2|, G_1>G_2 |T_1<T_2) \\
&=\frac{1}{2}P(G_2>G_1>|T_1-T_2| ) \\
&\quad + \frac{1}{2}P(G_1>|T_1-T_2|, |G_1-G_2|<|T_1-T_2|, G_1>G_2),
\ea
\edm
where $T_1,T_2$ (resp. $G_1,G_2$) are i.i.d. random variables following the distribution $F_T$ (resp. $F_G$).
and we use the symmetry $P(T_1<T_2)=P(T_2<T_1)=1/2$ in derivation of the above equation.
The results can be further simplified using symmetry:
\bdm
\ba
&P(G_2>G_1>|T_1-T_2|) = P(G_1>G_2>|T_1-T_2|) \\
&=\frac{1}{2}P(G_2>G_1>|T_1-T_2|)+\frac{1}{2}P(G_1>G_2>|T_1-T_2|) \\
&=\frac{1}{2}P(\min\{G_1,G_2\}>|T_1-T_2|,\;G_2>G_1) + \frac{1}{2}P(\min\{G_1,G_2\}>|T_1-T_2|,\;G_1>G_2) \\
&=\frac{1}{2}P(\min\{G_1,G_2\}>|T_1-T_2|),
\ea
\edm
and similarly,
\bdm
\ba
&P(G_1>|T_1-T_2|, |G_1-G_2|<|T_1-T_2|, G_1>G_2) \\
&=\frac{1}{2}P(\max\{G_1,G_2\}>|T_1-T_2|, |G_1-G_2|<|T_1-T_2|).
\ea
\edm
Therefore, the probability of the first event has the the form
\bdm
\ba
&P(T_u<T_v<T_u+G_u<T_v+G_v)=\frac{1}{4}P(\min\{G_1,G_2\}>|T_1-T_2| ) \\
&\quad + \frac{1}{4}P(\max\{G_1,G_2\}>|T_1-T_2|, |G_1-G_2|<|T_1-T_2|).
\ea
\edm

By symmetry again, the probability of the second event is equal to that of the first event.
The probability of the third event can be computed similarly as
\bdm
\ba
&P(T_u<T_v<T_v+G_v<T_u+G_u)=P(T_1<T_2, G_1-G_2>T_2-T_1) \\
&=P(G_1-G_2>T_2-T_1|T_1<T_2)\cdot P(T_1<T_2) \\
&=\frac{1}{2}P(G_1-G_2>|T_1-T_2||T_1<T_2) \\
&=\frac{1}{2}P(G_1-G_2>|T_1-T_2||T_1<T_2,G_2<G_1)\cdot P(G_2<G_1) \\
&=\frac{1}{4}P(|G_1-G_2|>|T_1-T_2||T_1<T_2,G_2<G_1) \\
&=\frac{1}{4}P(|G_1-G_2|>|T_1-T_2|).
\ea
\edm
Therefore, the expectation of $\Delta$ is given as 
\bdm
\ba
\E[\Delta]&=P(T_u<T_v<T_u+G_u<T_v+G_v)
 + P(T_v<T_u<T_v+G_v<T_u+G_u) \\
&\quad + 2 P(T_u<T_v<T_v+G_v<T_u+G_u) \\
&=2P(T_u<T_v<T_u+G_u<T_v+G_v) + 2 P(T_u<T_v<T_v+G_v<T_u+G_u) \\
&=\frac{1}{2}P(\min\{G_1,G_2\}>|T_1-T_2| ) + \frac{1}{2}P(\max\{G_1,G_2\}>|T_1-T_2|, |G_1-G_2|<|T_1-T_2|) \\
&\quad + \frac{1}{2}P(|G_1-G_2|>|T_1-T_2|).
\ea
\edm
The expectation of $|\cS_u|$ is then equal to
\bdm
\E[|\cS_u|] = 1 + \sum_{v\in\cV\setminus\{u\}}\E[\Delta_v]=1+(|\cV|-1)\mu,
\edm
where $\Delta_v$ is the increment of time buckets on $|\cS_u|$ induced by $\rtw_v$.
The tail bound \eqref{eqn:P(Su)} is a direct application of the Hoeffding's inequality.
\halmos
\endproof

\corUniform*
\proof{Proof.}
Define the following random variables:
\bdm
\Tdiff=|T_1-T_2|,\;\Gdiff=|G_1-G_2|,\;\Gmin=\min\{G_1,G_2\},\;\Gmax=\max\{G_1,G_2\}.
\edm
Let $f_{\Tdiff}(\cdot)$, $f_{\Gdiff}(\cdot)$ and $f_{\Gmin}(\cdot)$ be the probability density functions
of $\Tdiff$, $\Gdiff$ and $\Gmin$, respectively. 
Using the convolution formula, one can verify that
\bdm
\ba
&f_{\Tdiff}(x)=\frac{2(H-x)}{H^2} & x\in[0,H], \\
&f_{\Gdiff}(x)=f_{\Gmin}(x)=\frac{2(h-x)}{h^2} & x\in[0,h].
\ea
\edm
We compute each term in \eqref{eqn:mu}.\allowdisplaybreaks
\beq
\ba
&P(\Gmin>\Tdiff)=\int^h_0f_{\Tdiff}(x)\cdot\left(\int^h_x f_{\Gmin}(y)dy\right)dx \\
&=\int^h_0\frac{2(H-x)}{H^2}\cdot\left(\int^h_x\frac{2(h-y)}{h^2}dy\right)dx \\
&=\int^h_0\frac{2(H-x)}{H^2}\cdot\frac{(h-x)^2}{h^2}dx \\
&=\int^h_0\frac{2(H-h+z)z^2}{H^2h^2}dz \\
&=\int^h_0\frac{2(H-h)z^2}{H^2h^2}dz + \int^h_0\frac{2z^3}{H^2h^2}dz \\
&=\frac{2h}{3H}-\frac{h^2}{6H^2}.
\ea
\eeq
Since $\Gdiff$ and $\Gmin$ has the same density function, we have $P(\Gdiff>\Tdiff)=P(\Gmin>\Tdiff)$.
\beq
\ba
&P(\Gdiff<\Tdiff<\Gmax)=2P(\Gdiff<\Tdiff<\Gmax,G_1<G_2) \quad (\textrm{by symmetry})\\
&=2P(G_2-G_1<\Tdiff<G_2,G_1<G_2) \\
&=\int^h_0\int^h_{x_1}\frac{2}{h^2}\left(\int^{x_2}_{x_2-x_1}f_{\Tdiff}(y)dy\right)dx_1dx_2 \\
&=\int^h_0\int^h_{x_1}\frac{2}{h^2}\left(\int^{x_2}_{x_2-x_1}\frac{2(H-y)}{H^2}dy\right)dx_1dx_2 \\
&=\int^h_0\int^h_{x_1}\frac{2}{h^2}\left[\frac{(H-(x_2-x_1))^2}{H^2} - \frac{(H-x_2)^2}{H^2} \right]dx_1dx_2 \\
&=\int^h_0\left(\int^h_{x_1}\frac{2(H-(x_2-x_1))^2}{H^2h^2}dx_2\right)dx_1 - \int^h_0\left(\int^h_{x_1}\frac{2(H-x_2)^2}{H^2h^2}dx_2\right)dx_1\\
&=\int^h_0\left(\frac{2H}{3h^2}-\frac{2(H-h+x_1)^3}{3H^2h^2} \right)dx_1 - \int^h_0\left(\frac{2(H-x_1)^3}{3H^2h^2}-\frac{2(H-h)^3}{3H^2h^2} \right)dx_1\\
&=\frac{2H}{3h}-\frac{H^2}{6h^2}+\frac{(H-h)^4}{6H^2h^2} + \frac{(H-h)^4}{6H^2h^2}-\frac{H^2}{6h^2} + \frac{2(H-h)^3}{3H^2h} \\
&=\frac{2h}{3H}-\frac{h^2}{3H^2}.
\ea
\eeq
Finally, we obtain $\mu$ in this case as 
\bdm
\ba
\mu&=\frac{1}{2}P(\Gmin>\Tdiff)+\frac{1}{2}P(\Gdiff<\Tdiff<\Gmax)+\frac{1}{2}P(\Gdiff>\Tdiff) \\
&=P(\Gmin>\Tdiff)+\frac{1}{2}P(\Gdiff<\Tdiff<\Gmax) \\
&=\frac{h}{H}-\frac{h^2}{3H^2},
\ea
\edm
concluding the proof.
\halmos
\endproof

\section{Problem instance generation for numerical experiments}
\label{sec:prob-gen}
The 30 tree-FRCVP instances can be divided into
two major groups, each containing 15 instances. 
We set $\lambda=\infty$ (recall that $\lambda$ is the capacity of a platoon) 
for instances in the first group, and $\lambda=10$ in the second group.
The 15 instances within each group are generated corresponding to 15 combinations
of the network topology (three options: ArtifNet (Art), Chicago-cluster (CC) and Chicago-full (CF)) 
and the vehicle number $N$ (5 options: 100, 200, 300, 400 and 500). 
The 15 instances are labeled by the network type and $N$. 
For example, Art200 represents the instance of a 200-vehicle
system on the ArtifNet, and the labels of all instances are provided in the first two columns of
Table~\ref{tab:comp-perf}. 
However, the route graph defined by shortest paths of all vehicles is not necessarily a tree,
and hence some minor adjustments of routes are made to recast the graph to be a tree.
Examples of tree route graphs are given in Figure~\ref{fig:network} 
(b), (d) and (f) for a 100-vehicle system on the three networks, respectively. 

A tree-FRCVP problem instance consists of three parts: the road network topology,
the origin and destination nodes of each vehicle ($O_v$ and $D_v$), and the origin and destination
time of each vehicle ($T^O_v$ and $T^D_v$). 
\subsection{Road networks generation}
We create three road networks 
as the base for all problem instances. The first one is an artificial network that is 
generated using an attachment random graph model. The network 
generation process starts from a single node locating at $(0,0)$ and 
an arc is added to a randomly selected node (named as an active node)
with the degree less than four. 
This network is denoted as ArtifNet in the section, and the network structure is shown 
in Figure~\ref{fig:network}(a), which consists of 100 directed edges.

The second and third road networks are both based on the Chicago high-way
network in reality. The third road network referred as Chicago-full (CF)
is the Chicago high-way network itself, shown in Figure~\ref{fig:network}(e). 
The network has about 36K nodes and 48K edges. 
The second road network is a simplified version of the Chicago-full network by imposing
clustering on it. Specifically, we apply a $K$-mean clustering algorithm with $K=100$ to
the nodes of Chicago-full network with the distance of two nodes measured by the
$\ell_2$ norm of the two positions. So all nodes from the Chicago-full network are partitioned
into 100 clusters, and then for nodes in a same cluster, a representative node is created
at the geometric center of nodes in the cluster. The 100 representative nodes are connected
with newly generated edges to form a planar graph shown in Figure~\ref{fig:network}(c),
and it is referred as the Chicago-cluster (CC) network in the section.

\subsection{Vehicle number, origin and destination nodes}
The second part of creating a problem instance is to specify the number of vehicles
in the system and generate the origin and destination nodes of each vehicle.
We consider five options for the number $N$ of vehicles in a network: 100, 200, 300, 400 and 500. 
For a given network and vehicle number $N$, we select $N$ pairs of origin and destination nodes, 
each for a vehicle. Note that vehicles could share a same origin (resp., destination) node. 
For the ArtifNet network, the origin (resp., destination) nodes are randomly selected from the negative 
(resp., positive) side of horizontal axis. For the Chicago-cluster and Chicago-full networks,
the origin (resp., destination) nodes are randomly selected from the Northern (resp., Southeastern) 
part of the network. We set the route of a vehicle to be the shortest path from the origin 
to the destination in all the numerical instances investigated in this paper. 
Obviously there are many ways for the route selection. In particular, the routes can 
also be taken as an optimal solution of the coordinated vehicle platooning problem 
without time constraints, i.e., the solution of (RDP) in \cite{luo2020-route-then-sched-vehicle-platn}.  
An example of origin and destination nodes for a 100-vehicle system 
is given in Figure~\ref{fig:network}(b) for the ArtifNet network, Figure~\ref{fig:network}(d) for the
Chicago-cluster network, and Figure~\ref{fig:network}(f) for the Chicago-full network, respectively.
Restricting origin (resp., destination) nodes to a certain portion of the network may lead to 
large overlap of vehicle routes, which is in favor of coordinated platooning. 

\subsection{Origin and destination time}\label{sec:EC-origin-dest-time}
The last step is to generate origin and destination time (i.e., $T^O_v$ and $T^D_v$) 
for each vehicle. To simplify the implementation, we assume that the travel speed of 
every vehicle on every road link is a constant, and hence the traversing time of a road link
can be simply set to be the length of the road link. In order to generate reasonable origin 
and destination time for each vehicle, we first compute the mean value $\ol{L}$ of 
route length over all vehicles, i.e., $\ol{L}=(\sum^N_{v=1}L_v)/N$, where $L_v$ is
the route length of vehicle $v$. Based on two input parameters $\gamfull$,
and $\gamext$ ($\gamfull>\gamext$), we generate $T^O_v$ and $T^D_v$ 
independently for each vehicle by the following two steps: First, draw a sample
from the uniform distribution $\mtc{U}[0,\gamfull\ol{L}]$ and set it to be $T^O_v$,
and then we set $T^O_v\gets T^O_v+(1+\gamext)L_v$. The parameter $\gamfull$ controls
the full scope of the time horizon while $\gamext$ controls the portion of flexible time relative
to the travel time. For the instances generated for computational performance comparison 
discussed in Section~\ref{sec:comp-perf}, we choose $\gamfull=50$ and $\gamext=2$.
For the platoon capacity $\lambda$, we consider two options $\lambda\in\{10,\infty\}$,
where $\lambda=\infty$ means no restriction on the platoon size.

\end{document}